\input ./style/arxiv-general.cfg
\documentclass[MSNbibl,number,citesort,seceqn,dvips]{arxbj}
\makeatletter
   \@ifpackageloaded{graphicx}{}{\usepackage{graphicx}}
\makeatother
\usepackage{tikz,upgreek}
\aid{0}
\volume{21}
\issue{4}
\pubyear{2015}
\firstpage{2513}
\lastpage{2551}
\doi{10.3150/14-BEJ653} % kopijuoti is 'New paper accepted'
\docsubty{FLA}

\makeatletter

\newcommand{\lleft}{\left}
\newcommand{\rrvert}{\vert}
\newcommand{\rright}{\right}
\newcommand{\rrVert}{\Vert}
\newcommand{\llvert}{\vert}
\newcommand{\llVert}{\Vert}
\newcommand{\Max}{
\mathbin{
\mathchoice
{\buildMax{\displaystyle}}
{\buildMax{\textstyle}}
{\buildMax{\scriptstyle}}
{\buildMax{\scriptscriptstyle}}
}
}

\newcommand{\Kuch}{
\mathbin{
\mathchoice
{\buildKuch{\displaystyle}}
{\buildKuch{\textstyle}}
{\buildKuch{\scriptstyle}}
{\buildKuch{\scriptscriptstyle}}
}
}

\newcommand\buildKuch[1]{%
\begin{tikzpicture}[baseline=(X.base), inner sep=0, outer sep=0]
\node[draw,circle] (X) {$#1\alpha$};
\end{tikzpicture}%
}

\def\Max{\bigcirc\hspace{-2.94mm}\vee}
\def\Exp{\operatorname{Exp}}
\def\Ee{\mathbf{E}}

\newtheorem{theorem}{Theorem}[section]
\newtheorem{lemma}[theorem]{Lemma}
\newtheorem{prop}[theorem]{Proposition}
\newproclaim{definition}[theorem]{Definition}
\newremark{remark}[theorem]{Remark}
\newremark{ex}{Example}[section]
\newremark{example}{Example}
\newcommand{\eqref}[1]{(\ref{#1})}

\newcommand{\sign}{\operatorname{sign}}
\renewcommand{\pi}{\uppi}
\renewcommand{\emptyset}{\varnothing}
\renewcommand{\leqslant}{\leq}
\renewcommand{\geqslant}{\geq}

\def\sfrac#1#2{#1/#2}
\def\vfrac#1#2{(#1)/#2}

\def\sklfrac#1#2{(#1/#2)}

\renewcommand{\epsilon}{\varepsilon}

\makeatother

\begin{document}
\begin{frontmatter}

\title{L\'evy processes and stochastic integrals in the sense of
generalized convolutions}
\runtitle{Weak L\'evy processes}

\begin{aug}
%%%% inicialai - be tarpu
\author[A]{\inits{M.}\fnms{M.}~\snm{Borowiecka-Olszewska}\thanksref{A}\ead[label=e1]{m.borowiecka-olszewska@wmie.uz.zgora.pl}},
\author[B]{\inits{B.H.}\fnms{B.H.}~\snm{Jasiulis-Go{\l}dyn}\corref{}\thanksref{B}\ead[label=e2]{jasiulis@math.uni.wroc.pl}},
\author[C]{\inits{J.K.}\fnms{J.K.}~\snm{Misiewicz}\thanksref{C}\ead[label=e3]{J.Misiewicz@mini.pw.edu.pl}}
\and
\author[D]{\inits{J.}\fnms{J.}~\snm{Rosi\'{n}ski}\thanksref{D}\ead[label=e4]{rosinski@math.utk.edu}}
%\author{\inits{}\fnms{}~\snm{}}
%%\runauthor{} %% auto
%\dedicated{}
\address[A]{Faculty of Mathematics, Computer Science and Econometrics,
University of Zielona G{\'o}ra, ul. Prof. Z.~Szafrana 4A, 65-516
Zielona G{\'o}ra, Poland.\\ \printead{e1}}
\address[B]{Institute of Mathematics, University of Wroc{\l}aw, pl.
Grunwaldzki 2/4, 50-384 Wroc{\l}aw, Poland.\\ \printead{e2}}
\address[C]{Faculty of Mathematics and Information Science, Warsaw
University of Technology, ul. Koszykowa 75, 00-662 Warszawa, Poland.
\printead{e3}}
\address[D]{Department of Mathematics, 227 Ayres Hall, University of Tennessee,
Knoxville, TN 37996, USA. \printead{e4}}
\end{aug}

% HISTORY:
\received{\smonth{12} \syear{2013}}
\revised{\smonth{3} \syear{2014}}

% ABSTRACT
%
\begin{abstract}
In this paper, we present a comprehensive theory of generalized and
weak generalized convolutions, illustrate it by a large number of
examples, and discuss the related infinitely divisible distributions.
We consider L\'evy and additive process with respect to generalized and
weak generalized convolutions as certain Markov processes,
and then study stochastic integrals with respect to such processes.
We introduce the representability property of weak generalized convolutions.
Under this property and the related weak summability, a stochastic
integral with respect to random measures related to such convolutions
is constructed.
\end{abstract}

% KEYWORDS
% visi is mazosios raides ir pagal abecele
%
\begin{keyword}
\kwd{L\'evy process}
\kwd{scale mixture}
\kwd{stochastic integral}
\kwd{symmetric stable distribution}
\kwd{weakly stable distribution}
\end{keyword}
\end{frontmatter}

%s1 #&#
\section{Introduction}\label{sec1}

Motivated by the seminal work of Kingman \cite{King}, K. Urbanik
introduced and developed the theory of generalized convolutions in his
fundamental papers \cite{Urbanik64,Urbanik73,Urbanik84,Urbanik86}.
Roughly speaking, a generalized convolution is a binary associative
operation $\star$ on probability measures such that the convolution of
point-mass measures $\delta_x \star\delta_y$ can be a nondegenerate
probability measure, while the usual convolution gives $\delta_{x+y}$.
The study of weakly stable distributions, initiated by Kucharczak and
Urbanik (see \cite{KU,Urbanik76}) and followed by a series of
papers by Urbanik, Kucharczak, Panorska, and Vol'kovich (see, e.g.,
\cite{KU2,APan1,APan2,vol1,vol2,vol3}), provided
a new and rich class of weak generalized convolutions on $\mathbb
{R}_{+}$ (called also $\mathcal{B}$-generalized convolutions).
Misiewicz, Oleszkiewicz and Urbanik \cite{MOU} gave full
characterization of weakly stable distributions with nontrivial
discrete part and proved some uniqueness properties of weakly stable
distributions that will be used in this paper. For additional
information on generalized convolutions and weakly stable laws, see
\cite{Basia,KendallWalk,JasKula,weaklevy,MisJas,JarMis,Log,Mis2006,MM,KO,APan3}.

%%%%%%%%%%%%%%%%%%%%%%%%%%tu skonczylam przeklejac

In this paper, we present a comprehensive theory of generalized and
weak generalized convolutions and discuss the related classes of
infinitely divisible distributions. We construct L\'evy and additive
processes with respect to such convolutions. L\'evy process with
respect to generalized convolutions form interesting subclasses of
Markov processes, such as the class of Bessel processes in the case of
Kingman's convolution (see \cite{Thu}), but in general, they are heavy
tailed Markov processes (see Remark~\ref{h-tail}). Then we construct
stochastic integrals of deterministic functions associated with such
convolutions and the corresponding L\'evy processes.
We also introduce the weak summability property of generalized
convolutions. If a convolution admits the weak summability, then the
stochastic integration theory related to such convolutions becomes more
explicit and concrete.

This paper is organized as follows. In Section~\ref{sec2}, we give definitions
and properties of generalized and weak generalized convolutions that
will be used throughout this work. We also provide an extensive list of
examples. In Section~\ref{sec3}, we recall main results on infinite divisibility
with respect to generalized and weak generalized convolutions. This
information is crucial for further considerations.
The main result of Section~\ref{sec4} states that under minimal assumptions on
generalized convolutions an analog of processes with independent
increments can be constructed. We follow and extend an approach of N.
Van Thu \cite{Thu94}. In Section~\ref{sec5}, we consider stochastic integral
processes with respect to generalized convolutions.
Section~\ref{sec6} is devoted to the property of weak generalized summation. In
Section~\ref{sec7}, we construct ``independently scattered'' random measures based
on a weak generalized summation; these measures are used in Section~\ref{sec8}
to construct L\'evy and additive processes. Finally, in Section~\ref{sec9} we
define stochastic integrals of deterministic functions with respect to
such random measures and generalized convolutions.

Throughout this paper, the distribution of the random
element $X$ is denoted by $\mathcal{L}(X)$. If $\lambda= \mathcal
{L}(X)$ and $a \in\mathbb{R,}$ we denote the law of $aX$ by $T_a
\lambda$. ${\mathcal P}(\mathbb{E})$ denotes the family of all
probability measures on the Borel $\sigma$-algebra $\mathcal
{B}(\mathbb{E})$ of a Polish space $\mathbb{E}$. For short, we write
${\mathcal P}(\mathbb{R})={\mathcal P}$ and $\mathcal{P}(\mathbb
{R}_{+}) = \mathcal{P}_{+}$. The set of all symmetric probability
measures on $\mathbb{R}$ is denoted by ${\mathcal P}_s$. If $\lambda
\in{\mathcal P}$ and $\lambda= \mathcal{L}(\theta)$, then $|\lambda
| \in{\mathcal P}_{+}$ is defined by $|\lambda| = \mathcal
{L}(|\theta|)$. If $\mu= \mathcal{L}(X)$ and $\lambda= \mathcal
{L}(\theta)$ are such that $X$ and $\theta$ are independent, then by
$\mu\circ\lambda$ we denote the distribution of $X \theta$.

%s2 #&#
\section{Generalized convolutions} \label{sec2}

%s2.1 #&#
\subsection{Urbanik's generalized convolutions}\label{sec2.1}
Urbanik \cite{Urbanik64} introduced a \emph{generalized
convolution} as a binary, symmetric, associative and commutative
operation $\diamond$ on ${\mathcal{P}}_{+}$ having the following properties:
\begin{enumerate}[(iii)]
\item[(i)] $\lambda\diamond\delta_0 = \lambda$ for all $\lambda
\in\mathcal{P}_{+}$;
\item[(ii)] $ ( p \lambda_1 + (1-p)\lambda_2  )
\diamond\lambda= p  ( \lambda_1 \diamond
\lambda ) + (1-p)  ( \lambda_2 \diamond
\lambda )$ for each $p \in[0,1]$ and $\lambda, \lambda_1,
\lambda_2 \in\mathcal{P}_{+}$;
\item[(iii)] $T_a  ( \lambda_1 \diamond\lambda_2  )
=  (T_a \lambda_1 ) \diamond (
T_a \lambda_2 )$ for all $a\geqslant0$ and $\lambda_1, \lambda
_2 \in\mathcal{P}_{+}$;
\item[(iv)] if $\lambda_n \rightarrow\lambda$ and $\nu_n
\rightarrow\nu$, then $(\lambda_n \diamond\nu_n) \rightarrow
(\lambda\diamond\nu)$, where $\rightarrow$ denotes the weak convergence;
\item[(v)] there exists a sequence of positive numbers $(c_n)$ such
that $T_{c_n} \delta_1^{\diamond n}$ converges weakly to a measure
$\nu\neq\delta_0$ (here $\lambda^{\diamond n} = \lambda\diamond
\cdots\diamond\lambda$ denotes the generalized convolution of $n$
identical measures $\lambda$).
\end{enumerate}
The property (v) is important. It states that for the generalized
convolution a kind of limit theorem holds with a nontrivial limit
measure. Another important property, which follows from (ii) and (iv),
is that for every $\lambda_1, \lambda_2 \in\mathcal{P}_{+}$ and a
Borel set $A \subset\mathbb{R}_{+}$
%
%e2.1 #&#
\begin{equation}
\label{kernel1} \lambda_1 \diamond\lambda_2 (A) = \int
_0^{\infty} \int_0^{\infty
}
(\delta_x \diamond\delta_y) (A) \lambda_1
(\mathrm{d}x) \lambda_2(\mathrm{d}y)
\end{equation}
(see Lemma~\ref{l:kernel2} for the proof of a related equality).
In view of \eqref{kernel1}, in order to specify $\diamond$ we only
need to know $\delta_x \diamond\delta_y$ for all $x,y$. Actually, it
is enough to know $\delta_z \diamond\delta_1$ for all $z \in[0,1]$,
because $\delta_x \diamond\delta_y = T_x ( \delta_1 \diamond\delta
_{y/x})$ for any $x>y$.

\subsection*{Examples}
For details, see \cite{Bingham,CKS,JasKula,King,KU2,Urbanik64,Urbanik73,Urbanik76,Urbanik84,Urbanik86,Urbanik88,vol1}.
\setcounter{ex}{-1}

\begin{ex}\label{ex2.0}
 The classical convolution (\cite{Urbanik64,Urbanik88}) is evidently an example of generalized
convolution. It will be denoted as usual by $\ast$:
\[
\delta_a \ast\delta_b = \delta_{a+b}.
\]
\end{ex}

\begin{ex}\label{ex2.1}
 \emph{Symmetric} generalized convolution (\cite{Urbanik64,Urbanik88}) on $\mathcal{P}_{+}$ is defined by
\[
\delta_a \ast_{s} \delta_b =
\tfrac{1}{2} \delta_{|a-b|} + \tfrac
{1}{2} \delta_{a+b}.
\]
The name \emph{symmetric} comes from the fact that this convolution
can be easily extended to a generalized convolution on $\mathcal{P}$
taking values in the set of symmetric measures $\mathcal{P}_{s}$:
\[
\delta_a \ast_{s} \delta_b =
\tfrac{1}{4} \delta_{a-b} + \tfrac
{1}{4} \delta_{-a+b}
+ \tfrac{1}{4} \delta_{-a-b} + \tfrac{1}{4}
\delta_{a+b}.
\]
\end{ex}

\begin{ex}\label{ex2.2}
 In a similar way another generalized convolution
(called by Urbanik $(\alpha,1)$-convolution in \cite{Urbanik64,Urbanik76}) can be defined for every $\alpha>0$ by means of
\[
\delta_a \ast_{s, \alpha} \delta_b =
\tfrac{1}{2} \delta _{|a^{\alpha}-b^{\alpha}|^{1/{\alpha}}} + \tfrac{1}{2} \delta
_{(a^{\alpha}+b^{\alpha})^{1/{\alpha}}}.
\]
\end{ex}

\begin{ex}\label{ex2.3}
 For every $p \in(0,\infty]$, the formula
\[
\delta_a \ast_{p} \delta_b =
\delta_c,\qquad a,b \geqslant0, \qquad c = \bigl\| (a,b)\bigr\|_p =
\bigl(a^p + b^p \bigr)^{1/p}
\]
defines a generalized convolution $\ast_{p}$ ($p$-stable convolution)
on $\mathcal{P}_{+}$. For details, see \cite{Urbanik64,Urbanik86}.
\end{ex}

\begin{ex}\label{ex2.4}
 The Kendall convolution $\diamond_{\alpha}$ on
$\mathcal{P}_{+}$, $\alpha> 0$, is defined (\cite{JasKula}) by
\[
\delta_x \diamond_{\alpha} \delta_1 =
x^{\alpha} \pi_{2\alpha} + \bigl(1-x^{\alpha} \bigr)
\delta_1, \qquad x\in[0,1],
\]
where $\pi_{2\alpha}$ is a Pareto measure with density $g_{2\alpha
}(x) = 2\alpha x^{-2\alpha-1} \mathbf{1}_{[1,\infty)}(x)$.
\end{ex}

\begin{ex}\label{ex2.5}
 The Kingman convolution $\otimes_{\omega_s}$ on
$\mathcal{P}_{+}$, $s> - \frac{1}{2}$, is defined in \cite{King} by
\[
\delta_a \otimes_{\omega_s} \delta_b =
\mathcal{L} \bigl(\sqrt{ a^2 + b^2 + 2ab
\theta_s } \bigr),
\]
where $\theta_s$ is absolutely continuous with the density function
\[
f_s (x)= \frac{\Gamma(s+1)}{\sqrt{\pi} \Gamma(s + \sfrac{1}{2})} \bigl( 1 - x^2
\bigr)_{+}^{s - \sfrac{1}{2}}.
\]
If $n:= 2(s+1) \in\mathbb{N}$, $n >1$, the variable $\theta_s$ can
be interpreted as one dimensional projection of the random vector
$\mathbf{U} = (U_1,\dots,U_n)$ having uniform distribution $\omega
_n$ on the unit sphere $S_{n-1} \subset\mathbb{R}^n$. If $n=1$ and $s
= - \frac{1}{2}$, then $\theta_s$ has the discrete distribution
$\frac{1}{2} \delta_{-1} + \frac{1}{2} \delta_1$.
\end{ex}

\begin{ex}\label{ex2.6}
 $\infty$-convolution (\cite{KU2,Urbanik88}) is
defined by
\[
\delta_a \mathbin{\Max}\,\delta_b =
\delta_{\max\{a,b\}}.
\]
\end{ex}

\begin{ex}\label{ex2.7}
 A combination of Kingman convolution and $(\alpha
,1)$ convolution, called by Urbanik $(\alpha, \beta)$-convolution in
\cite{Urbanik64}, for $0<\alpha<\infty, 0 < \beta< \infty$, is
defined for $a,b>0$ as
\[
\delta_a \otimes_{\alpha,\beta}\delta_b =
\mathcal{L} \bigl( \bigl(a^{2\alpha} + b^{2\alpha} + 2 a^{\alpha}
b^{\alpha} \theta \bigr)^{1/{2\alpha}} \bigr),
\]
where $\theta= \theta_{\vfrac{\beta-2}{2}}$ is a random variable
with the density function
\[
f_{(\beta-2)/2} (x) = \frac{\Gamma({\beta}/2)}{\sqrt{\pi} \Gamma
({(\beta-1)}/2)} \bigl( 1 - x^2
\bigr)_{+}^{(\beta-3)/2}.
\]
\end{ex}

\begin{ex}\label{ex2.8}
 A kind of generalization of Kendall convolution
called the Kucharczak--Urbanik convolution (\cite{Bingham}) was
obtained by the following definition for $\alpha>0$ and $s \in[0,1]$
\begin{eqnarray*}
\delta_s \diamond_{\alpha,n} \delta_1 (
\mathrm{d}x) &=& \bigl( 1 - s^{\alpha
} \bigr)_{+}^n
\delta_1 (\mathrm{d}x)
\\
&&{} + \frac{\alpha(n+1) s^{\alpha(n+1)}}{x^{2 \alpha n +1}} \sum_{k=1}^n {n
\choose{k}} {n\choose{k-1}} \frac{ (x^{\alpha} -
s^{\alpha} )^{k-1}  ( x^{\alpha} - 1
)_{+}^{n-k}}{s^{\alpha k} } \,\mathrm{d}x.
\end{eqnarray*}
\end{ex}

\begin{ex}\label{ex2.9}
The Kucharczak convolution $\Kuch$, $\alpha\in
(0,1)$, is defined in \cite{Urbanik88} by
\[
\delta_a \Kuch \delta_b (\mathrm{d}x) =
\frac{a^{\alpha}
b^{\alpha} \sin(\pi\alpha)(2x - a -b)}{\pi(x-a-b)^{\alpha}
(x-a)^{\alpha} (x-b)^{\alpha}} \mathbf{1}_{ ( (a^{\alpha} +
b^{\alpha})^{1/{\alpha}}, \infty )}(x) \,\mathrm{d}x.
\]
\end{ex}

\begin{ex}\label{ex2.10}
The Vol'kovich convolution $\vartriangle
_{1,\beta}$ for $0<\beta< \frac{1}{2}$ (see \cite{vol1}) is given by
\[
\delta_a \vartriangle_{1,\beta} \delta_b (
\mathrm{d}x) = \frac{2 a^{2\beta
} b^{2\beta}}{B(\beta, \sklfrac{1}{2} - \beta)} \bigl[ \bigl( x^2 -
(a-b)^2 \bigr)_{+} \bigl( (a+b)^2 -
x^2 \bigr)_{+} \bigr]^{-\beta
- \sfrac{1}{2}} \,\mathrm{d}x.
\]
\end{ex}

\begin{ex}\label{ex2.11}
 In \cite{KU2} for $\alpha\in(0,1)$, the
authors considered the following measure:
\[
\mu= \bigl(2 - 2^{-\alpha} \bigr) \sum_{n=0}^{\infty}
2^{-1-n(\alpha+1)} T_{2^n} ( \pi_{\alpha} ),
\]
where $\pi_{\alpha}$ is the Pareto distribution with the density
$\alpha x^{-\alpha- 1} \mathbf{1}_{[1, \infty)}(x)$.
They proved that for every pair $a,b>0$ there exists a unique
probability measure $\varrho(a,b)\in\mathcal{P}_{+}$ fulfilling the equality
\[
T_a (\mu) \mathbin{\Max }\,T_b (\mu) = \mu\circ
\varrho(a,b).
\]
Setting $\delta_a \nabla_{\alpha} \delta_b := \varrho(a,b)$ they
obtained a generalized convolution. In a similar way, many other
generalized convolutions can be constructed on the basis of known
convolutions (see, e.g., \cite{JasKula}).
\end{ex}

\begin{ex}\label{ex2.12}
 We say that the distribution $\mu$ on $\mathbb
{R}^n$ is $\ell_1$-symmetric (sometimes the name $\ell
_1$-pseudo-isotropic is used here) if the characteristic function of
$\mu$ has the following form
\[
\widehat{\mu}(\xi) = \varphi \bigl( \|\xi\|_1 \bigr),
\]
for some function $\varphi$, where $\| \xi\|_1 = |\xi_1| +
\cdots+ |\xi_n|$. This means that the random vector $\mathbf{X}$ is
$\ell_1$-symmetric ($\ell_1$-pseudo-isotropic) if for every
$\xi\in\mathbb{R}^n$ the following equation holds
\[
\langle \xi, \mathbf{X}\rangle = \sum_{k=1}^n
\xi_k X_k \stackrel{\mathrm{d}} {=} \|\xi
\|_1 \cdot X_1.
\]
In 1983,  Cambanis,  Keener and  Simons \cite{CKS} described
the set of extreme points of the family of $\ell_1$-symmetric
distributions on
$\mathbb{R}^n$. They proved that the random vector $\mathbf{X}$ is
$\ell_1$-pseudo-isotropic iff there exists a nonnegative random
variable $\Theta$ such that
%
%e2.2 #&#
\begin{equation}
\label{1-iso} \mathbf{X} \stackrel{\mathrm{d}} {=} \biggl( \frac{U_1}{\sqrt{D_1}},
\dots, \frac{U_n}{\sqrt{D_n}} \biggr) \cdot\Theta=: \mathbf {V}\cdot\Theta,
\end{equation}
where $\mathbf{U}^n = (U_1,\dots,U_n)$ has uniform
distribution on the unit sphere in $\mathbb{R}^n$, $\mathbf{D} =
(D_1, \dots, D_n)$ has Dirichlet distribution with
parameters $(\frac{1}{2}, \dots, \frac{1}{2})$, $\mathbf{U}^n$,
$\mathbf{D}$ and $\Theta$ are independent. This means that the set of
extreme points for the set of $\ell_1$-pseudo-isotropic distributions
on $\mathbb{R}^n$ is equal to
\[
\bigl\{ T_a \mathcal{L}(\mathbf{V}) \colon a \geqslant0 \bigr\}.
\]

Let $\varphi(\| \xi\|_1)$ be the characteristic function of $\mathbf
{V}$, that is, $\varphi(\| \xi\|_1) = \mathbf{E} \mathrm{e}^{\mathrm{i}\langle \xi, \mathbf
{V}\rangle}$. Then the characteristic function of $a \mathbf{V} + b \mathbf
{V'}$, where $\mathbf{V'}$ is an independent copy of a $\mathbf{V}$,
is of the form
\[
\Phi \bigl(\| \xi\|_1 \bigr) = \varphi \bigl(a \| \xi\|_1 \bigr) \varphi
\bigl(b\| \xi\|_1 \bigr),
\]
thus it also depends only on $\|\xi\|_1$. By \eqref{1-iso}, there
exists a random variable $\Theta= \Theta(a,b)$ independent of
$\mathbf{V}$ such that
\[
a \mathbf{V} + b \mathbf{V'} \stackrel{\mathrm{d}} {=} \mathbf{V}
\Theta.
\]
Now we obtain a generalized convolution $\nabla_{\ell_1}$ setting
\[
\delta_a \nabla_{\ell_1} \delta_b =
\mathcal{L} \bigl( \Theta (a,b) \bigr).
\]
Unfortunately, an explicit formula for $\mathcal{L} ( \Theta
(a,b)  )$ is unknown.
\end{ex}

%
%re2.1 #&#
\begin{remark}
By Schoenberg's classical result (see \cite{Schoenberg}), we have that
a random vector $\mathbf{X}$ on $\mathbb{R}^n$ is $\ell
_2$-pseudo-isotropic ($\ell_2$-symmetric, rotationally invariant) iff
$\mathbf{X} \stackrel{\mathrm{d}}{=} \mathbf{U} \sqrt{\Theta}$ for some
nonnegative variable $\Theta$ independent of $\mathbf{U}$. This leads
to the family of Kingman's convolutions in special cases $n=2(s+1) \in
\mathbb{N}$. The characterization \eqref{1-iso} proven in \cite{CKS}
gives a general form for $\ell_1$-pseudo-isotropic and leads to the
generalized convolution $\nabla_{\ell_1}$. In both cases the
distributions of the extreme points of $\ell_i$-pseudo-isotropic
measures, $i=1,2$, that is, $\mathcal{L}(\mathbf{U})$ and $\mathcal
{L}(\mathbf{V})$ are weakly stable. A full characterization of $\ell
_{\alpha}$-symmetric distributions for $\alpha\notin\{1,2\}$ is
unknown. All we know is that only $\alpha\leqslant2$ can be
considered here.
\end{remark}

A pair $(\mathcal{P}_+,\diamond)$ is called a \textit{generalized
convolution algebra}. A continuous mapping $h\dvtx  \mathcal{P}_+ \to
\mathbb{R}$ is called a \textit{homomorphism} of $(\mathcal
{P}_+,\diamond)$ if
\begin{itemize}[$\bullet$]

\item[$\bullet$] $\forall a\in[0,1]\ \forall \lambda_1,\lambda_2
\in\mathcal{P}_+\ h( a \lambda_1 + (1-a) \lambda_2 ) = a h(\lambda
_1) + (1-a) h(\lambda_2)$,
\item[$\bullet$] $\forall \lambda_1, \lambda_2 \in\mathcal{P}_+
\ h( \lambda_1 \diamond\lambda_2 ) = h(\lambda_1) h(\lambda_2)$.
\end{itemize}
Obviously, $h(\cdot)\equiv0$ and $h(\cdot)\equiv1$ are the trivial
homomorphisms. A generalized convolution algebra $(\mathcal
{P}_+,\diamond)$ is said to be \emph{regular} if it admits a
nontrivial homomorphism.
%

%de2.2 #&#
\begin{definition}
We say that a nontrivial generalized convolution algebra $(\mathcal
{P}_+,\diamond)$ admits a characteristic function if there exists
one-to-one correspondence $\lambda\leftrightarrow\Phi_{\lambda}$
between probability measures $\lambda\in\mathcal{P}_{+}$ and real
valued functions $\Phi_{\lambda}$ on $[0,\infty)$ such that for
$\lambda, \nu\in\mathcal{P}_{+}$
\begin{enumerate}[3.]
\item[1.] $\Phi_{p \lambda + q \nu} = p \Phi_\lambda + q \Phi
_{\nu}$ for $p,q \geqslant0$, $p+q=1$;
\item[2.] $\Phi_{\lambda\diamond\nu} = \Phi_{\lambda} \cdot\Phi
_{\nu}$;
\item[3.] $\Phi_{T_a \lambda}(t) = \Phi_{\lambda} (at)$;
\item[4.] the uniform convergence of $\Phi_{\lambda_n}$ on every
bounded interval is equivalent to the weak convergence of $\lambda_n$.
\end{enumerate}
The function $\Phi_{\lambda}$ is called the characteristic function
of the probability measure $\lambda$ in the algebra $(\mathcal
{P}_+,\diamond)$ or $\diamond$-generalized characteristic function of
$\lambda$.
\end{definition}

It can be shown (see \cite{Urbanik84}) that $\Phi$ is uniquely
determined up
to a scale parameter.

The $\diamond$-generalized characteristic function in generalized
convolution algebra plays the same role as the classical Laplace or
Fourier transform for convolutions defined by addition of independent
random elements. The following fact is crucial for further
investigations, see \cite{Urbanik64} for the proof.

%
%pr2.3 #&#
\begin{prop}\label{cf}
A nontrivial generalized convolution algebra $(\mathcal{P}_+,\diamond
)$ admits a characteristic function $\Phi$ if and only if it is
regular. In this case
\[
\Phi_{\lambda} (t) = h (T_t \lambda ),\qquad t\geqslant0,
\lambda \in\mathcal{P}_{+},
\]
where $h$ is the nontrivial homomorphism of $(\mathcal{P}_+,\diamond
)$. Moreover, the map $\lambda\mapsto\Phi_{\lambda}$ is an integral
transform:
\[
\Phi_{\lambda}(t) = \int_0^{\infty} \Omega(tx)
\lambda(\mathrm{d}x),
\]
where $\Omega(t) := h(\delta_t)$. $\Omega$ is called the kernel of
the $\diamond$-generalized characteristic function $\Phi$.
\end{prop}

It can be shown that for each nontrivial homomorphism $h$ on a regular
algebra $(\mathcal{P}_+,\diamond)$ there exists an open neighborhood
of zero $U$ such that
\[
\forall x \in U \setminus\{0\},\qquad 0< \bigl|h(\delta_x)\bigr|<1.
\]
This property implies that the $\diamond$-generalized characteristic
function $\Phi_{\lambda}(\cdot)$ of the measure $\lambda\in
\mathcal{P}_+$ has a very useful property: if $\Phi_{\lambda}(t_n) =
1$ for some $t_n \searrow0$, then $\lambda= \delta_0$. One can find
more about generalized convolutions in \cite{JasKula,KU,KU2,Kill,Urbanik64,Urbanik73,Urbanik76,Urbanik84,Urbanik86,vol1,vol2,vol3}.

%s2.2 #&#
\subsection{Weak generalized convolutions}\label{sec2.2}
Weak generalized convolutions were studied in \cite{Basia,weaklevy,Mis2006,MM,MOU,Urbanik76,vol1}. They are derived from the concept
of weakly stable probability measures.

%
%de2.4 #&#
\begin{definition}
The distribution $\mu$ of a random vector $\mathbf{X}$, taking values
in a separable Banach space $\mathbb{E}$, is weakly stable if for
every $a,b \in\mathbb{R}$ there exists a random variable $\theta$
independent of $\mathbf{X}$ such that
{\renewcommand{\theequation}{$\ast$}
\begin{eqnarray}
a \mathbf{X}_1 + b \mathbf{X}_2 \stackrel{d} {=} \theta
\mathbf{X},
\end{eqnarray}}
where $\mathbf{X}_1, \mathbf{X}_2$ are independent copies of $\mathbf
{X}$ and $\stackrel{d}{=}$ denotes equality in distribution.
\end{definition}

If the condition ($\ast$) holds only for nonnegative constants $a,b$,
then we say that $\mathbf{X}$ is $\mathbb{R}_{+}$-\textit{weakly
stable}. It was shown in \cite{MOU} that if a weakly stable measure
$\mu$ has an atom, then either $\mu=\delta_0$ or $\mu= \frac{1}{2}
\delta_{\mathbf{a}} + \frac{1}{2} \delta_{-\mathbf{a}}$ for some
$\mathbf{a} \in\mathbb{E}$. In both cases we shall call such
measures trivial.

It was proved in \cite{MOU} that the condition ($\ast$) is equivalent
to the following:
{\renewcommand{\theequation}{$\ast\ast$}
\begin{eqnarray}
\forall \theta_1, \theta_2\  \exists \theta\qquad
\theta_1 \mathbf{X}_1 + \theta_2
\mathbf{X}_2 \stackrel{\mathrm{d}} {=} \theta\mathbf{X},
\end{eqnarray}}%
where $\theta_1, \theta_2$ are random variables such that $\theta_1,
\theta_2, \mathbf{X}_1, \mathbf{X}_2$ are independent and $\theta$
is independent of $\mathbf{X}$. Then (\ref{astast}) can be written in
the language of distributions in the following way:
\[
(\mu\circ\lambda_1) \ast(\mu\circ\lambda_2) = \mu\circ
\lambda,
\]
where $\mathcal{L}(\theta_i) = \lambda_i$, $i=1,2$, $\mathcal
{L}(\theta) = \lambda$. If the measure $\mu$ is nonsymmetric, then
$\lambda$ is uniquely determined from $\mu\circ\lambda$, but when
$\mu$ is symmetric, then only the measure $|\lambda| = \mathcal
{L}(|\theta|)$ (equivalently, $\frac{1}{2} (\lambda+ T_{-1}\lambda
)$) is uniquely determined (see \cite{MOU}).

Having a weakly stable random vector $\mathbf{X}$ with distribution
$\mu$, we are able to define a weak generalized convolution:

%
%de2.5 #&#
\begin{definition}
Let $\mu\in{\mathcal P}(\mathbb{E})$ be a nontrivial weakly stable
measure, and let $\lambda_1, \lambda_2\in{\mathcal P}$. If
\[
(\mu\circ\lambda_1 ) \ast (\mu\circ\lambda_2 ) = \mu
\circ\lambda,
\]
then the weak generalized convolution (also called $\mu$-weak
generalized convolution) of the measures $\lambda_1, \lambda_2$ with
respect to the measure $\mu$
(notation $\lambda_1 \otimes_{\mu}\lambda_2$) is defined as follows
\[
\lambda_1 \otimes_{\mu} \lambda_2 = \lleft\{
\begin{array} {l@{\qquad}l} \lambda& \mbox{if $\mu$ is not symmetric;}
\\
|\lambda| & \mbox{if $\mu$ is symmetric.} \end{array} %
\rright.
\]
Sometimes it is more convenient to define $\lambda_1 \otimes_{\mu}
\lambda_2 = \frac{1}{2} (\lambda+ T_{-1}\lambda)$, when $\mu$ is
symmetric. The pair $( \mathcal{P}, \otimes_{\mu})$ is called a weak
generalized convolution algebra.
\end{definition}

The following lemma describes basic properties of weak generalized convolution.

%
%le2.6 #&#
\begin{lemma}\label{prop-weak}
If the weakly stable measure $\mu\in{\mathcal P}(\mathbb{E})$ is not
trivial, then for all $\lambda, \lambda_1, \lambda_2, \lambda_3 \in
\mathcal{P}$
\renewcommand\theenumi{(\arabic{enumi})}
\renewcommand\labelenumi{\theenumi}
\begin{enumerate}[(6)]
\item\label{p1} $\lambda_1 \otimes_{\mu} \lambda_2$ is uniquely
determined;
\item\label{p2} $\lambda_1 \otimes_{\mu}\lambda_2 = \lambda_2
\otimes_{\mu}\lambda_1$;
\item\label{p3} $ (\lambda_1\otimes_{\mu}\lambda_2 )
\otimes_{\mu} \lambda_3 = \lambda_1 \otimes_{\mu}
 ( \lambda_2\otimes_{\mu}\lambda_3)$;
\item\label{p4} $\lambda\otimes_{\mu}\delta_0 = \lambda$
($\lambda
\otimes_{\mu}\delta_0 = |\lambda|$ if $\mu$ is
symmetric); %
\item\label{p5} $ ( p \lambda_1 + (1-p)\lambda_2  )
\otimes_{\mu}\lambda= p  ( \lambda_1 \otimes_{\mu}
\lambda ) + (1-p)  ( \lambda_2 \otimes_{\mu}
\lambda )$ for each $p \in[0,1]$;
\item\label{p6} $T_a  ( \lambda_1 \otimes_{\mu}\lambda_2
 )
=  (T_a \lambda_1 ) \otimes_{\mu} (
T_a \lambda_2 )$;
\item\label{p7} if $\lambda_n \rightarrow\lambda$ and $\nu_n
\rightarrow\nu$, then $\lambda_n \otimes_{\mu} \nu_n \rightarrow
\lambda\otimes_{\mu} \nu$.
\end{enumerate}
\end{lemma}

\begin{pf} Property~\ref{p1} follows from Theorems 3 and 4 in
\cite{MOU}. Properties \ref{p2}--\ref{p6} are simple
consequences of the definition and the uniqueness property \ref{p1}.
To see \ref{p7} it is enough to notice that for independent random
sequences $Y, Y_1, Y_2, \dots$ and $Z, Z_1, Z_2, \dots$ the following
implications hold
\[
Y_n \stackrel{\mathrm{d}} {\rightarrow} Y,\qquad Z_n
\stackrel{\mathrm{d}} { \rightarrow} Z \quad \Rightarrow\quad \lleft\{ %
\begin{array} {l} Y_n \cdot Z_n \stackrel{\mathrm{d}} {
\rightarrow} Y \cdot Z;
\\
Y_n + Z_n \stackrel{\mathrm{d}} {\rightarrow} Y + Z,
\end{array} %
\rright.
\]
where $\stackrel{\mathrm{d}}{\rightarrow}$ denotes convergence of
distributions, and then use the uniqueness \ref{p1}.
\end{pf}

%
%le2.7 #&#
\begin{lemma} \label{l:kernel2}
If $\mu$ is not symmetric, then for every $\lambda_1, \lambda_2 \in
\mathcal{P}$ and $A \in\mathcal{B}(\mathbb{R})$ we have
\setcounter{equation}{2}%
%e2.3 #&#
\begin{equation}
\label{kernel2} \lambda_1 \otimes_{\mu}
\lambda_2 (A) = \int_{\mathbb{R}^2} (\delta_x
\otimes_{\mu} \delta_y) (A) \lambda_1 (
\mathrm{d}x) \lambda_2(\mathrm{d}y).
\end{equation}
If $\mu$ is symmetric, then this equality holds with $\mathbb{R}$
replaced by $\mathbb{R}_+$ and $\lambda_1, \lambda_2 \in\mathcal{P}_+$.
\end{lemma}

\begin{pf} If $\lambda_1, \lambda_2$ have finite supports,
$\lambda_1 = \sum_{i=1}^m p_i \delta_{x_i}$, $\lambda_2 = \sum_{j=1}^n q_j \delta_{y_j}$, then by \ref{p2} and \ref{p5},
$\lambda_1 \otimes_{\mu} \lambda_2 = \sum_{i,j} p_i q_j \delta
_{x_i} \otimes_{\mu} \delta_{y_j}$. Hence for a bounded continuous
function $f$ on $\mathbb{R,}$ we have
%
%e2.4 #&#
\begin{eqnarray}
\label{kern-f} \int_{\mathbb{R}} f(z) (\lambda_1
\otimes_{\mu} \lambda_2) (\mathrm{d}z) & =& \sum
_{i,j} p_i q_j \int
_{\mathbb{R}} f(z) (\delta_{x_i} \otimes _{\mu}
\delta_{y_j}) (\mathrm{d}z)
\nonumber
\\[-8pt]
\\[-8pt]
& =& \int_{\mathbb{R}^2} \int_{\mathbb{R}} f(z) (
\delta_x \otimes_{\mu
} \delta_y) (\mathrm{d}z)
\lambda_1 (\mathrm{d}x) \lambda_2(\mathrm{d}y).
\nonumber
\end{eqnarray}
Let $\lambda_1, \lambda_2 \in\mathcal{P}$ be arbitrary. Choose
$\lambda_{i,n} \in\mathcal{P}$ with finite supports such that
$\lambda_{i,n} \to\lambda_i$ as $n \to\infty$, $i=1,2$.
We have for any bounded continuous function $f$ on $\mathbb{R}$
\begin{eqnarray*}
&&\int_{\mathbb{R}^2} \int_{\mathbb{R}} f(z) (
\delta_x \otimes _{\mu} \delta_y) (
\mathrm{d}z) \lambda_1 (\mathrm{d}x) \lambda_2(
\mathrm{d}y)
\\
&&\quad = \lim_{n\to\infty} \int_{\mathbb{R}^2} \int
_{\mathbb{R}} f(z) (\delta_x \otimes_{\mu}
\delta_y) (\mathrm{d}z) \lambda_{1,n} (\mathrm{d}x)
\lambda _{2,n}(\mathrm{d}y)
\\
&&\quad = \lim_{n\to\infty} \int_{\mathbb{R}} f(z) (
\lambda_{1,n} \otimes_{\mu} \lambda_{2,n}) (
\mathrm{d}z) = \int_{\mathbb{R}} f(z) (\lambda _1
\otimes_{\mu} \lambda_2) (\mathrm{d}z).
\end{eqnarray*}
The first equality holds because the map $(x,y) \mapsto\int f(z)
(\delta_x \otimes_{\mu} \delta_y)(\mathrm{d}z)$ is continuous by \ref{p7}
and bounded, the second one follows from \eqref{kern-f}, and the third
uses \ref{p7}. We have shown
\[
\int_{\mathbb{R}} f(z) (\lambda_1
\otimes_{\mu} \lambda_2) (\mathrm{d}z) = \int
_{\mathbb{R}^2} \int_{\mathbb{R}} f(z) (
\delta_x \otimes_{\mu
} \delta_y) (\mathrm{d}z)
\lambda_1 (\mathrm{d}x) \lambda_2(\mathrm{d}y)
\]
for any bounded continuous function $f$. By a standard monotone class
argument, we deduce that this equality holds for any $f= \mathbf
{1}_A$, $A \in\mathcal{B}(\mathbb{R})$, which gives \eqref
{kernel2}. The proof in the symmetric case of $\mu$ is similar.
\end{pf}

Notice that for a weak generalized convolution the condition (v) of the
Urbanik definition of generalized convolution does not have to be
satisfied. In \cite{Basia}, we can find a wide description
of properties of the generalized convolutions on $\mathbb{R}$ without
property (v). However it was shown in \cite{MOU} that if the measure
$\mu$ has a finite weak moment of order $\varepsilon> 0$, then there
exists a measure $\lambda$ such that $\mu\circ\lambda$ is symmetric
$\alpha$-stable for some (and then for
every) $\alpha\leq\min\{\varepsilon, 2\}$. This means that $T_{c_n}
\lambda^{\otimes_{\mu} n} = \lambda$ for a properly
chosen sequence $(c_n)$, and the property (v) holds if we replace
$\delta_1$ by $\lambda$.

The weak generalized convolution is always regular with
\[
\Omega(t)= h(\delta_t) := \widehat{\mu}(t) = \int
_{\mathbb{R}} \mathrm{e}^{\mathrm{i}tx} \mu(\mathrm{d}x),
\]
see Proposition~\ref{cf}.

\subsection*{Examples}
{\renewcommand{\theexample}{2.1a}
\begin{example}\label{ex2.1a}
Let $\theta$ be a random variable with
distribution $\lambda_0 = \frac{1}{2}\delta_1 + \frac{1}{2} \delta
_{-1}$ and let $\theta'$ be its independent copy. It is easy to check
that for all $a,b \geqslant0$, $a\neq b$
\[
a \theta+ b \theta' \equiv \bigl| a \theta+ b \theta' \bigr|
\cdot\frac{a \theta+ b \theta'}{|a \theta+ b \theta'|},
\]
where the two factors on the right are independent and
\[
\frac{a \theta+ b \theta'}{|a \theta+ b \theta'|} \stackrel{\mathrm{d}} {=} \theta.
\]
This shows that $\theta$ is weakly stable. Moreover, since
\[
\mathcal{L} \bigl(\bigl|a \theta+ b \theta'\bigr| \bigr) = \tfrac{1}{2}
\delta_{|a-b|} + \tfrac{1}{2} \delta_{a+b},
\]
we have that the symmetric generalized convolution is a weak
generalized convolution and $\ast_s = \otimes_{\lambda_0}$.
\end{example}}

{\renewcommand{\theexample}{2.3a}
\begin{example}\label{ex2.3a}
 Not for all $p>0$, but for $p \in(0,2]$ the
convolution $\ast_p$ can be extended to a weak generalized convolution
on $\mathcal{P}$ taking values in $\mathcal{P}_{+}$ defined by
$\gamma_p$-symmetric $p$-stable measure which is weakly stable since
\[
a \Gamma_p + b \Gamma_p' \equiv\bigl\|(a,b)
\bigr\|_p \Gamma_p'',\qquad \mbox
{where } \Gamma_p'' := \frac{a}{\|(a,b)\|_p}
\Gamma_p + \frac{b}{\|
(a,b)\|_p} \Gamma_p,
\]
where $\Gamma_p, \Gamma_p'$ are independent with the distribution
$\gamma_p$. Evidently, the first equality holds everywhere and, by the
basic properties of stable variables, $\Gamma_p''$ also has the
distribution~$\gamma_p$.
\end{example}}

{\renewcommand{\theexample}{2.4a}
\begin{example}\label{ex2.4a}
Not for all $\alpha>0$, but for $\alpha\in
(0,1]$ the Kendall convolution $\diamond_{\alpha}$ can be extended to
a weak generalized convolution on $\mathcal{P}$ taking values in
$\mathcal{P}_{s}$ defined by the measure $\mu_{\alpha}$ with the
characteristic function $\widehat{\mu_{\alpha}} (t) = ( 1 -
|t|^{\alpha})_{+}$.
\end{example}}

{\renewcommand{\theexample}{2.5a}
\begin{example}\label{ex2.5a}
 For $2(s+1) \in\mathbb{N}$ the Kingman
convolution has the natural interpretation as a weak generalized
convolution with respect to the weakly stable uniform distribution on
the unit sphere $S_{2s+1} \subset\mathbb{R}^{2(s+1)}$. More
precisely:

Let $\mathbf{U}_n$, $n \geqslant2$, denotes the random vector with
the uniform distribution $\omega_n$ on the unit sphere $S_{n-1}
\subset\mathbf{R}^n$. It is known that if $\mathbf{U},\mathbf{U}'$
are independent copies of $\mathbf{U}_n$, then for each $a,b \in
\mathbb{R}$, $ab\neq0$, the random variables
\[
\bigl\llVert a\mathbf{U} + b \mathbf{U}' \bigr\rrVert
_2 \quad \mbox{and}\quad \frac{
a\mathbf{U} + b \mathbf{U}'}{\llVert  a\mathbf{U} + b \mathbf
{U}'\rrVert _2}
\]
are independent and the second one has the distribution $\omega_n$. Since
\[
a\mathbf{U} + b \mathbf{U}' \equiv \bigl\llVert a\mathbf{U} + b
\mathbf {U}' \bigr\rrVert _2 \frac{ a\mathbf{U} + b \mathbf{U}'}{\llVert
a\mathbf{U} + b \mathbf{U}'\rrVert _2} \qquad
\mbox{a.e.}
\]
this implies that $\omega_n$ is weakly stable and it defines the
weakly stable convolution $\otimes_{\omega_n}$ on $\mathcal{P}$ in
the following way
\[
\delta_a \otimes_{\omega_n} \delta_b =
\mathcal{L} \bigl( \bigl\llVert a\mathbf{U} + b \mathbf{U}' \bigr
\rrVert _2 \bigr).
\]
For $2s=-1$, we simply have
\[
\delta_a \otimes_{\omega_1} \delta_b =
\tfrac{1}{2} \delta_{|a-b|} + \tfrac{1}{2} \delta_{a+ b},
\]
which is $\ast_{s}$ convolution considered in Example~\ref{ex2.1}.
\end{example}}

{\renewcommand{\theexample}{2.12a}
\begin{example}\label{ex2.12a}
 By the result of Cambanis, Keener and Simons
\cite{CKS}, the distribution of $\mathbf{V}$ is weakly stable, and by
our construction
\[
\nabla_{\ell_1} = \otimes_{\mathcal{L}(\mathbf{V})}.
\]
\end{example}}

%s3 #&#
\section{Infinite divisibility with respect to generalized convolutions}\label{sec3}

%s3.1 #&#
\subsection{Infinite divisibility (decomposability) of measures on
$\mathbb{R}_{+}$}\label{sec3.1}

It is natural to consider infinitely divisible measures with respect to
generalized convolutions. Following Urbanik \cite{Urbanik64},
sometimes we will call such measures \emph{infinitely decomposable}.

%
%de3.1 #&#
\begin{definition}
A measure $\lambda\in\mathcal{P}_+$ is said to be infinitely
divisible with respect to the generalized convolution $\diamond$
($\diamond$-infinitely decomposable) in the algebra $(\mathcal
{P}_{+}, \diamond)$ if for every $n \in\mathbb{N}$ there exists a
probability measure $\lambda_n \in\mathcal{P}_+$ such that $\lambda
= \lambda_n^{\diamond n}$.
\end{definition}

The proof of the following proposition can be found in \cite{Urbanik64}.

%
%pr3.2 #&#
\begin{prop}\label{prop1}
Let $\lambda\in\mathcal{P}_{+}$ be $\diamond$-infinitely divisible.
There exists a collection of measures $\lambda^{\diamond r}$,
$r\geqslant0$ such that
\begin{enumerate}[(iii)]
\item[(i)] $\lambda^{\diamond0} = \delta_0$, $\lambda^{\diamond1}
= \lambda$;
\item[(ii)] $\lambda^{\diamond r} \diamond\lambda^{\diamond s} =
\lambda^{\diamond(r+s)}$, $r,s \geqslant0$;
\item[(iii)] $\lambda^{\diamond r_n} \rightarrow\delta_0$ if $r_n
\searrow0$.
\end{enumerate}
\end{prop}

Similarly as in the classical theory, one of the most important
examples of $\diamond$-infinitely divisible distribution is given by
\[
\Exp_{\diamond} (a \lambda) \stackrel{\mathrm{def}} {=}
\mathrm{e}^{-a} \sum_{k=0}^{\infty}
\frac{a^k}{k!} \lambda^{\diamond k},
\]
where $\lambda\in\mathcal{P}_+$ and $a>0$. The measure $\Exp
_{\diamond} (a \lambda)$ is called a \textit{generalized compound
Poisson measure} or $\diamond$-\textit{compound Poisson measure}. If
$\lambda= \delta_1$, then it is called a \textit{generalized Poisson
measure} or $\diamond$-\textit{Poisson measure}. To see that $\Exp
_{\diamond} (a \lambda)$ is infinitely divisible with respect to
$\diamond$ it is sufficient to observe that
\[
\biggl(\Exp_{\diamond} \biggl( \frac{a}{n} \lambda \biggr)
\biggr)^{\diamond n}=\Exp_{\diamond} (a \lambda).
\]

Another important example of a $\diamond$-infinitely divisible
distribution gives the following:

%
%de3.3 #&#
\begin{definition}
Let $\lambda\in\mathcal{P}_+$. We say that $\lambda$ is stable in
the generalized convolution algebra $(\mathcal{P}_+, \diamond)$ if
the following condition holds:
\[
\forall a,b \geqslant0\ \exists c\geqslant0 \qquad T_a \lambda
\diamond T_b \lambda= T_c \lambda.
\]
\end{definition}

%
%re3.4 #&#
\begin{remark}
A measure $\lambda$ is stable in the generalized convolution algebra
$(\mathcal{P}_+, \diamond)$ (or simply $\diamond$-\textit{stable})
if and only if there exists a sequence of positive numbers $(c_n)$ and
$\eta\in\mathcal{P}_+$ such that
\[
T_{c_n} \eta^{\diamond n} \to\lambda.
\]
For details of the proof see Theorem~14 in \cite{Urbanik64}.
\end{remark}

In the formulation of the analog of the L\'evy--Khintchine formula for a
$\diamond$-infinitely divisible distribution we need the
characteristic exponent $\varkappa(\diamond)$ for the generalized
convolution $\diamond$ defined in the following theorem of Urbanik
\cite{Urbanik86}:

%
%th3.5 #&#
\begin{theorem}
For every generalized convolution $\diamond$ on $\mathcal{P}_{+}$
there exists a constant $\varkappa(\diamond) \in(0, \infty]$ such
that for every $p \in(0, \varkappa(\diamond)]$ there exists a
measure $\sigma_p \in\mathcal{P}_{+}$ with the $\diamond
$-generalized characteristic function
\[
\Phi_{\sigma_p}(t) = \lleft\{ %
\begin{array} {l@{\qquad}l}
\mathrm{e}^{-t^p} & \mbox{if } p < \infty;
\\
\mathbf{1}_{[0,1]}(t) & \mbox{if } p=\infty. \end{array} %
\rright.
\]
Moreover, the set of all $\diamond$-stable measures coincides with the set
\[
\bigl\{ T_a (\sigma_p) \colon a > 0, 0 < p \leqslant
\varkappa (\diamond) \bigr\}.
\]
\end{theorem}

In particular we have that
\[
\mathrm{e}^{-t^{\varkappa}} = \int_0^{\infty}
\Omega(ts) \sigma_{\varkappa}(\mathrm{d}s).
\]
Let $\gamma_{p+}$, $p \in(0,1)$, be the completely skewed to the
right stable measure $\gamma_{p+}\sim S_p(\sigma, 1,0)$ with $\sigma
^p = 2^p \cos\frac{p \pi}{2}$ and the Laplace transform $\mathrm{e}^{-2^p
t^p}$ with the notation $S_p(\sigma, \beta, \mu)$ as in the
representation 1.1.6 in \cite{ST} and let $\gamma_{p+}$ be the
distribution of $\theta_p$. Since
\[
\int_0^{\infty} \mathrm{e}^{{-t^{\varkappa}s}/2}
\gamma_{p+} (\mathrm{d}s) = \mathrm{e}^{-
t^{p \varkappa}},
\]
we see that for $s<\varkappa$ the measure $\sigma_s =
T_{2^{1/{\varkappa}}} \sigma_{\varkappa} \circ\mathcal{L}(\theta
_p^{1/{\varkappa}})$ for $p = \frac{s}{\varkappa}$ is absolutely
continuous with respect to the Lebesgue measure. In many cases also the
measure $\sigma_{\varkappa}$ is absolutely continuous with respect to
the Lebesgue measure.
We denote by $f_s$ the density
function for the standard $s$-stable measure with respect to the
generalized convolution $\diamond$.

It was proven (see Theorem~7 in \cite{Urbanik64}) by Urbanik that
the characteristic exponent does not depend on the choice of
nontrivial homomorphism and consequently on the choice of the
$\diamond$-generalized characteristic function.

The examples given below illustrate the material of this section.
Examples \ref{ex3.0} and \ref{ex3.3} belong to the classical theory of stable
distributions. Formulas for densities in Examples \ref{ex3.4} and \ref{ex3.11} are new.
Detailed calculations related to Examples \ref{ex3.5}--\ref{ex3.11} can mostly be found
in \cite{Urbanik88}. The formula for the generalized characteristic
function in Example~\ref{ex3.12} is new.

\subsection*{Examples}

{\renewcommand{\theexample}{3.0}
\begin{example}\label{ex3.0}
 As a nontrivial homomorphism in the case of
usual convolution on $\mathcal{P}_{+}$ we can simply take $h(\lambda)
= \int_0^{\infty} \mathrm{e}^{-x} \lambda(\mathrm{d}x)$, that is, the kernel of the
transform can be given by $\Omega(t) = h(T_t \delta_1) =
\mathrm{e}^{-t}\mathbf{1}_{[0,\infty)}(t)$. Moreover $\varkappa(\ast) = 1$,
$\sigma_1 = \delta_1$ and $\sigma_p = \gamma_{p+}$. In particular,
for $p=\frac{1}{2}$ we have $\sigma= 1$ and the density of $\gamma
_{\sfrac{1}{2}+}$ can be written in terms of elementary functions, namely
\[
\gamma_{\sfrac{1}{2}+} (\mathrm{d}x) = \frac{1}{\sqrt{2 \pi}} x^{- \sfrac
{3}{2}} \exp
\biggl\{ - \frac{1}{2x} \biggr\} \,\mathrm{d}x.
\]
It has been shown in \cite{Zolotarev}, that
\[
T_{\sfrac{1}{2}} \gamma_{\sfrac{1}{3}+} (\mathrm{d}x) = \frac{1}{3\pi}
x^{-\sfrac{3}{2}} K_{\sfrac{1}{3}} \biggl( \sqrt{ \frac{4}{27x}} \biggr) \,
\mathrm{d}x,
\]
where $K_{1/3}$ is the MacDonald function and
\[
T_{\sfrac{1}{2}} \gamma_{\sfrac{2}{3}+} (\mathrm{d}x) = \frac{1}{x \sqrt{3\pi
}}
W_{\sfrac{1}{2}, \sfrac{1}{6}} \biggl( \frac{4}{27x^2} \biggr) \exp \biggl\{ -
\frac{2}{27 x^2} \biggr\} \,\mathrm{d}x,
\]
where $W_{p,q}$ is the Whittaker function.
\end{example}}

{\renewcommand{\theexample}{3.3}
\begin{example}\label{ex3.3}
 For the generalized convolution $\ast_p$ on
$\mathcal{P}_{+}$ we have $\Omega(t) = \mathrm{e}^{-t^p}$, $\varkappa(\ast
_p) = p$ and $\sigma_s = \gamma_{\frac{s}{p}+}$ for $s<p$.
\end{example}}

{\renewcommand{\theexample}{3.4}
\begin{example}\label{ex3.4}
 For the generalized Kendall convolution $\diamond
_{\alpha}$ on $\mathcal{P}_{+}$ we have $\Omega(t) = (1-t^{\alpha
})_{+}$, $\varkappa(\diamond_{\alpha}) = \alpha$ and for $p \in
(0,\alpha]$
\[
f_p (x) = p x^{-p-1} \biggl( 1 - \frac{p}{\alpha} +
\frac{p}{\alpha
} x^{-p} \biggr) \mathrm{e}^{-x^{-p}}
\mathbf{1}_{(0,\infty)}(x).
\]
For the same convolution considered as an operation on $\mathcal
{P}_{s}$ the functions $\Omega$ and $f_p$ shall be symmetrized.
\end{example}}

{\renewcommand{\theexample}{3.5}
\begin{example}\label{ex3.5} For the generalized Kingman convolution, we have
\[
\Omega(t) = \Gamma ( s+1 ) \biggl( \frac{t}{2} \biggr)^{s}
J_{s} (t),
\]
where $J_{r}$ is the Bessel function, $\varkappa(\otimes_{\omega_s})
= 2$,
\[
f_1(x) = \frac{\Gamma(s+\sfrac{3}{2})}{\sqrt{\pi} \Gamma(s+1)} \frac{x^s \mathbf{1}_{(0, \infty)}(x)}{(1+x)^{s+ \sfrac{3}{2}}},\qquad
f_2(x) = \frac{\mathbf{1}_{(0, \infty)}(x)}{2^{s+1} \Gamma(s+1)} x^s \mathrm{e}^{ - {{x}/2}},
\]
and for $0<p<2$
\[
f_p(x) = \frac{x^s}{2^{s+1} \Gamma(s+1)} \int_0^{\infty}
y^{-s-1} \exp \biggl\{ - \frac{x}{2y} \biggr\}
\gamma_{\sfrac{p}{2}+}(\mathrm{d}y).
\]
\end{example}}

{\renewcommand{\theexample}{3.6}
\begin{example}\label{ex3.6}
For the $\infty$-convolution
\[
\Omega(t) = \mathbf{1}_{[0,1]}(t)
\]
and $\varkappa(\mathbin{\Max}\, ) =\infty$, $\sigma_{\infty} = \delta
_1$ and
\[
f_p(x) = p x^{-p-1} \exp \bigl\{ - x^{-p} \bigr\}
\mathbf{1}_{(0,\infty)}(x)
\]
is the Weibull--Gnedenko distribution. It has been proven by Urbanik
\cite{Urbanik86} that $\varkappa(\diamond) = \infty$ if and only if
$\diamond= \mathbin{\Max}\,$.
\end{example}}

{\renewcommand{\theexample}{3.9}
\begin{example}\label{ex3.9} For the Kucharczak convolution $\Kuch$, $\alpha
\in(0,1)$,
\[
\Omega(t) = \Gamma(\alpha)^{-1} \Gamma(\alpha, t),
\]
where $\Gamma(\alpha,t)$ is the incomplete Gamma function, $\varkappa
( \Kuch) = \alpha$ and
\[
\sigma_p \bigl([0,x )\bigr) = x^{1-\alpha} \int
_0^x (x-y)^{\alpha-1} \gamma_{p+}(
\mathrm{d}y).
\]
\end{example}}

{\renewcommand{\theexample}{3.10}
\begin{example}\label{ex3.10}
 For the Vol'kovich convolution with $0 < \beta<
\frac{1}{2}$ we have
\[
\Omega(t) = \frac{2^{1-\beta} t^{\beta}}{\Gamma(\beta)} K_{\beta}(t),
\]
where $K_{\beta}$ is the MacDonald function and $\varkappa(\triangle
_{1,\beta}) = 2 \beta$.
\end{example}}

{\renewcommand{\theexample}{3.11}
\begin{example}\label{ex3.11}
 For the generalized convolution $\nabla_{\alpha
}$, $\alpha\in(0,1)$, under $\infty$ -convolution we have
\[
\Omega(t) = \bigl( 1 - 2^{(1+\alpha)[\log_2t]} - \bigl(2 - 2^{-\alpha} \bigr)
\bigl(1- 2^{[\log_2t]} \bigr) t^{\alpha} \bigr) \mathbf{1}_{[0,1]}(t),
\]
where the square brackets denote the integer part and $\varkappa
(\nabla_{\alpha}) = \alpha$. Moreover
\[
\sigma_p \bigl([0,x )\bigr) = \frac{2^{1+\alpha}}{2^{1+\alpha}-1} \biggl(1+
\frac{p}{\alpha x^p} \biggr) \mathrm{e}^{ - x^{-p}} - \frac{1}{2^{1+\alpha
}-1}
\biggl(1+ \frac{p 2^p}{\alpha x^p} \biggr) \mathrm{e}^{ - 2^p x^{-p}}.
\]
\end{example}}

{\renewcommand{\theexample}{3.12}
\begin{example}\label{ex3.12}
 For the Cambanis, Keener and Simons convolution,
we have
\[
\Omega(t) = \mathbf{E} \mathrm{e}^{\mathrm{i}t V_1} = \frac{\Gamma(\sfrac
{n}{2})}{\sqrt{\pi}
\Gamma(\vfrac{n-1}{2})} \int
_1^{\infty} \Omega_n
\bigl(ur^2 \bigr) u^{-{n/2}} (u-1)^{{(n-3)}/2} \,\mathrm{d}u,
\]
where $\Omega_n(r^2)$ is the characteristic function of the first
coordinate of the vector $\mathbf{U}^n$ and  $\varkappa(\nabla_{\ell
_1}) = 1$. The measure $\sigma_p$, $p\leqslant1$, in this case is
such that $\mathcal{L}(\mathbf{V}) \circ\sigma_p = \gamma_p$, for
$\gamma_p$ being the symmetric $p$-stable measure (abbreviation: S$p
$S measure).
\end{example}}

The following theorem (see Theorem~13 in \cite{Urbanik64}) gives the
L\'evy--Khintchine formula for $\diamond$-generalized characteristic
function for a $\diamond$-infinitely divisible distribution.

%
%th3.6 #&#
\begin{theorem}\label{thm:Urbanik}
Let $( \mathcal{P}_+, \diamond)$ be a regular generalized convolution
algebra. A function $\Phi\dvtx \mathbb{R}_+ \rightarrow\mathbb{R}$
is a $\diamond$-generalized characteristic function of a $\diamond
$-infinitely divisible measure iff it has the following representation
\[
\Phi(t) = \exp \biggl\{ - A t^{\varkappa(\diamond)} + \int_{0}^{\infty}
\frac{ \Omega(tx) - 1 }{ \upsilon(x) } m (\mathrm{d}x) \biggr\},
\]
where $m$ is a finite Borel measure on $[0,\infty)$,
\[
\upsilon(x) = \lleft\{ %
\begin{array} {l@{\qquad}l} 1 - \Omega(x) &
\mbox{if } 0 \leqslant x \leqslant x_0,
\\
1 - \Omega(x_0) & \mbox{if } x \geqslant x_0 \end{array}
\rright.
\]
and $x_0 >0$ is such that $\Omega(x) <1$ whenever $0< x \leqslant x_0$.
\end{theorem}

%s3.2 #&#
\subsection{Weak infinite divisibility}\label{sec3.2}

It is known that if a weakly stable measure $\mu$ is symmetric and
such that $\int_{\mathbb{E}}| \langle\xi, \mathbf{x}\rangle |^{\varepsilon}
\mu(\mathrm{d}\mathbf{x}) < \infty$ for some $\varepsilon>0$ and all
continuous linear functionals $\xi\in\mathbb{E}^{\ast}$, then the
weak generalized convolution $\otimes_{\mu}$ is a generalized
convolution in the Urbanik sense (i.e., $\otimes_{\mu}$ has property
$(\mathrm{v})$), see \cite{weaklevy}. Consequently, the infinite
divisibility with respect to such convolutions on $\mathcal{P}_{+}$
was already described in the previous subsection. Adding the
information about weakly stable measures $\mu$, that generate such
convolutions, will make this description more detailed and concrete.

%
%de3.7 #&#
\begin{definition}
Let $\mu\in\mathcal{P}(\mathbb{E})$ be a weakly stable measure. We
say that the measure $\lambda$ is $\mu$-weakly
infinitely divisible if for every $n \in\mathbb{N}$ there exists a
probability measure $\lambda_n$ such that
\[
\lambda= \lambda_n^{\otimes_{\mu} n} \equiv\lambda_n^n
\stackrel {def} {=} \lambda_n \otimes_{\mu} \cdots
\otimes_{\mu} \lambda_n, \qquad \mbox{($n$-times),}
\]
where (for the uniqueness) $\lambda, \lambda_n \in{\mathcal P}_{+}$ if
$\mu$ is $\mathbb{R}_{+}$-weakly stable or if $\mu$ is symmetric, and
$\lambda, \lambda_n \in{\mathcal P}$ if $\mu$ is weakly stable
nonsymmetric.
\end{definition}

Notice that if $\lambda$ is ${\mu}$-weakly infinitely divisible, then
$\mu\circ\lambda$ is infinitely divisible in the classical sense.
This information can be of some help in investigations, however we
shall remember that the opposite implication does not hold. There are
measures $\lambda$ and weakly stable measures $\mu$ such that $\mu
\circ\lambda$ is infinitely divisible and $\lambda$ is not ${\mu
}$-weakly infinitely divisible. Counterexamples are known even for $\mu
$ symmetric Gaussian and symmetric stable measures $\mu$ (see Example~2 in \cite{weaklevy}). Special properties of infinitely divisible
sub-stable distributions are discussed in \cite{Log,Sato,ST}.

It was proven in \cite{weaklevy} that for every nontrivial weakly
stable measure $\mu$ and $\mu$-weakly infinitely divisible measure
$\lambda$ there exists a family of measures $\{\lambda^r: r \geq0\}$
such that
\begin{enumerate}[(3)]
\item[(1)] $\lambda^0 = \delta_0$, $\lambda^1 = \lambda$;

\item[(2)] $\lambda^r \otimes_{\mu} \lambda^s = \lambda^{r+s}$, $r,s \geq
0$;

\item[(3)] $\lambda^r \rightarrow\delta_0$ if $r\rightarrow0$.
\end{enumerate}

The $\mu$-\emph{weak compound Poisson measure} for the $\mu$-weak
generalized convolution is defined exactly in the same way (see \cite{weaklevy}) as the compound Poisson measure for generalized convolution:
\[
\Exp_{\otimes_{\mu}} (a \lambda) \stackrel{\mathrm{def}} {=}
\mathrm{e}^{-a} \sum_{k=0}^{\infty}
\frac{a^k}{k!} \lambda^{\otimes_{\mu} k},
\]
where $\lambda\in\mathcal{P}$ and $a>0$. Sometimes this measure is
called $\mu$-\emph{weak generalized exponent} of the measure $a
\lambda$. If $\lambda= \delta_1$, then it is called a \emph{$\mu
$-weak Poisson measure}. In the case of $\mu$-weak generalized
convolution the following additional interesting property holds:
\[
\mu\circ\Exp_{\otimes_{\mu}}(a\lambda) = \exp \bigl(a(\mu\circ \lambda) \bigr),
\]
that is, every $\mu$-weak compound Poisson measure is a factor of some
compound Poisson measure. In some cases, we get the explicit formulas
for the generalized Poisson distribution.

\subsection*{Examples}
{\renewcommand{\theexample}{3.3a}
\begin{example}\label{ex3.3a}
 Let $\mu= \gamma_p$, $p\in(0,2]$ be symmetric
$p$-stable distribution on $\mathbb{R}$ with the characteristic
function $\mathrm{e}^{-A|r|^p}$, $A>0$. Then the $\mu$-weak Poisson measure is
purely discrete with the distribution
\[
\Exp_{\otimes_{\mu}}(c \delta_1) = \mathrm{e}^{-c} \sum
_{k=0}^{\infty} \frac{c^k}{k!}
\delta_{k^{1/p}}.
\]
To see this, it is enough to notice that if $X_{1},\ldots,X_{k}$ are
independent random variables with distribution $\gamma_p$, then
$X_{1}+\cdots+X_{k}\stackrel{ \mathrm{d}}{=}k^{1/p}X_{1}$, thus
\[
\delta_1^{\otimes_{\gamma_p} k}=\delta_{k^{1/p}}.
\]
\end{example}}

{\renewcommand{\theexample}{3.4a$_s$}
\begin{example}\label{ex3.4as}
 Consider the Kendall weak generalized
convolution $\diamond_{\alpha}$ on $\mathcal{P}_s$ with respect to
the weakly stable measure $\mu_{\alpha}$ with the characteristic
function $\widehat{\mu_{\alpha}}(t)=  ( 1-|t|^{\alpha}
)_+$, $\alpha\in(0,1]$.
It was shown in \cite{KendallWalk} that
\[
\bigl( 1-|t|^{\alpha} \bigr)_+^{k}=\int_{\mathbb{R}}
\bigl( 1-|ts|^{\alpha} \bigr)_+\lambda_{k}(\mathrm{d}s),
\]
where $\lambda_0= \frac{1}{2} \delta_1 + \frac{1}{2} \delta_{-1}$
and for $k\geq1$ we have
\[
\lambda_{k}(\mathrm{d}s)=\frac{\alpha k (k-1 )}{2} \bigl( 1-
|s|^{-\alpha} \bigr)^{k-2}|s|^{-(2\alpha+1)} \mathbf{1}_{(1,\infty)}
\bigl(|s|\bigr)\,\mathrm{d}s.
\]
This means that $\delta_1^{\diamond_{\alpha}k}=\lambda_{k}$ for
$k\geqslant1$, thus $\mu_{\alpha} \circ\lambda_k = \mu_{\alpha
}^{\ast k}$, and the $\mu_{\alpha}$-weak generalized exponent of $c
\delta_1$ can be calculated as
\begin{eqnarray*}
\Exp_{\diamond_{\alpha}} (c \delta_1 ) (\mathrm{d}s)& =&
\mathrm{e}^{-c} \delta_0(\mathrm{d}s) +
\mathrm{e}^{-c} c \lambda_0 (\mathrm{d}s)
\\
&&{}+\mathrm{e}^{-c}\sum_{k=2}^{\infty}
\frac{c^k}{k!}\frac{\alpha k
(k-1 )}{2} \bigl( 1- |s|^{-\alpha}
\bigr)^
{k-2}|s|^{- (2\alpha+1 )}\mathbf{1}_{(1,\infty)}\bigl(|s|\bigr)\,
\mathrm{d}s
\\
&=& \mathrm{e}^{-c} ( \delta_0 + c \lambda_0 )
(\mathrm{d}s) + \frac{\alpha
c^2}{2|s|^{ (2\alpha+1 )}} \mathrm{e}^{- c|s|^{-\alpha}}
\mathbf{1}_{(1,\infty)}\bigl(|s|\bigr) \,\mathrm{d}s.
\end{eqnarray*}
\end{example}}

{\renewcommand{\theexample}{3.4a$_+$}
\begin{example}\label{ex3.4aplus}
 Consider the same Kendall weak
generalized convolution as an operator on $\mathcal{P}_{+}$. Then,
similarly as before for $\mathcal{P}_s$ case, we obtain that $\exp(c
\mu_{\alpha}) = \mu_{\alpha}\circ\Exp_{\diamond_{\alpha}}(c
\delta_1)$, where
\[
\Exp_{\diamond_{\alpha}}(c \delta_1) (\mathrm{d}u) =
\mathrm{e}^{-c} \delta_0 (\mathrm{d}u) + c
\mathrm{e}^{-c} \delta_1(\mathrm{d}u) + \frac{c^2 \alpha}{u^{2\alpha+1}}
\mathrm{e}^{-cu^{-\alpha}} \mathbf{1}_{(1,\infty)}(u) \,\mathrm{d}u.
\]
\end{example}}

{\renewcommand{\theexample}{3.5a}
\begin{example}\label{ex3.5a}
 For the technical reasons we consider here the
special case of the Kingman weak generalized convolution $\otimes
_{\omega_3} \dvtx \mathcal{P}_s \rightarrow\mathcal{P}_s$. Since
the generalized convolutions defined by $\omega_3$ and by its
one-dimensional projection $\omega_{3,1}$ are the same and $\omega
_{3,1}(\mathrm{d}u) = \frac{1}{2} \mathbf{1}_{[-1,1]}(u) \,\mathrm{d}u$, the calculations
are simpler than in the general case.

For any $c>0$ we need to calculate $\lambda= \Exp_{\otimes_{\omega
_{3,1}}}(\frac{1}{2}c \delta_1 + \frac{1}{2}c \delta_{-1})$ because
in $\mathcal{P}_s$ the role of $\delta_1$ is played by the measure
$\lambda_0 = \frac{1}{2} \delta_1 + \frac{1}{2} \delta_{-1}$.
Since $\widehat{\omega_{3,1}}(r) = \frac{\sin{r}}{r}$ and $\omega
_{3,1} \circ\lambda= \exp( c \omega_{3,1} )$ then
\[
\widehat{\omega_{3,1}\circ\lambda} (r) = \mathrm{e}^{ -c (1-\sfrac{\sin
r}{r} )}.
\]
On the other hand, we can write
\[
\widehat{\omega_{3,1}\circ\lambda} (r) = \int_{\mathbb
{R}}
\widehat{\lambda}(rs)\omega_{3,1}(\mathrm{d}s)= \int
_{-1}^1 \frac{1}{2}\widehat{\lambda}(rs) \,
\mathrm{d}s=\frac{1}{2r}\int_{-r}^r
\widehat{\lambda}(s) \,\mathrm{d}s.
\]
Thus,
\[
\int_{-r}^r \widehat{\lambda}(s) \,
\mathrm{d}s=2r \mathrm{e}^{ -c (1-\sfrac{\sin
r}{r} )}.
\]
From the last equation it follows that
\[
\widehat{\lambda}(r)=\frac{\mathrm{d}}{\mathrm{d}r} \bigl( r \mathrm{e}^{ -c (1-\sfrac
{\sin r}{r} )}
\bigr) =\mathrm{e}^{-c (1-\sfrac{\sin r}{r} )} \biggl( 1-{c} \frac{\sin
r}{r} + c \cos r
\biggr).
\]
This implies that
\[
\Exp_{\otimes_{\omega_{3,1}}}(c \lambda_0) = \exp (c \omega _{3,1}
) \ast ( \delta_0- c \omega_{3,1} + c \lambda_0
).
\]
\end{example}}

These examples show that the $\mu$-weak Poisson measure does not need
to be discrete, although it is a linear combination of $\mu$-weak
generalized convolutions of the Dirac measure $\delta_1$.

%
%de3.8 #&#
\begin{definition}
Let $\mu\in\mathcal{P}$ be a nontrivial weakly stable measure. A
measure $\lambda\in\mathcal{P} \setminus\{\delta_0\}$ is ${\mu
}$-weakly stable if there exists a sequence of positive numbers $(c_n)$
and a measure $\nu\in\mathcal{P}$ such that
\[
T_{c_n} \nu^{\otimes_{\mu} n} \rightarrow\lambda.
\]
\end{definition}

We denote by $\mathcal{S}(\mu)$ the set of all ${\mu}$-weakly stable
measures. Let
\[
\mathcal{S}_p (\mu) = \bigl\{ \lambda\in\mathcal{P} \setminus\{
\delta_0\} \colon T_a \lambda\otimes_{\mu}
T_b \lambda= T_{g_p(a,b)} \lambda \bigr\},
\]
where $g_p(a,b)= (|a|^p + |b|^p )^{1/p}$. Measures $\lambda$ in
$\mathcal{S}_p (\mu)$ will be
referred to as $\mu$-\textit{weakly} $p$-\textit{stable}.
For every symmetric weakly stable measure $\mu$ there exists a
parameter $\varkappa= \varkappa(\mu)$ called the \textit
{characteristic exponent}, such that
\[
\varkappa(\mu) =\sup \bigl\{ p \in(0,2] \colon\mathcal{S}_p(\mu)
\neq\emptyset \bigr\}.
\]
In our convention the supremum over the empty set equals zero. The
parameter $\varkappa$ is related to the symmetric $p$-stable measure
$\gamma_p$ in the usual sense. Note that $\varkappa(\mu) \leqslant
2$ for every weakly stable measure $\mu$ while the corresponding
characteristic exponent $\varkappa(\diamond)$ of the Urbanik type
generalized convolution $\diamond$ can take any value from the
positive half-line including infinity.
It was proven in \cite{weaklevy} that $\varkappa(\mu)$ has the
following characterization.

%
%th3.9 #&#
\begin{theorem}
For every weakly stable distribution $\mu$ and $\mathcal{M}(\mu) =
 \{ \mu\circ\lambda\colon\lambda\in\mathcal{P} \}$ we have
\begin{eqnarray*}
\varkappa(\mu) & = & \sup \biggl\{ p \in[0,2] \colon\int_{\mathbb{R}}|x|^p
\mu(\mathrm{d}x) < \infty \biggr\}
\\
& = & \sup \bigl\{ p \in[0,2] \colon\gamma_p \in \mathcal{M}(\mu)
\bigr\}.
\end{eqnarray*}
\end{theorem}

The next theorem gives us the analogue of the L\'{e}vy--Khintchine
representation for infinitely divisible distributions in the sense of
weak generalized convolution on $\mathcal{P}_s$. Here $\mathbb{R}_0 =
\mathbb{R} \setminus\{0\}$.

%
%th3.10 #&#
\begin{theorem}
Assume that $\mu$ is a nontrivial symmetric weakly stable measure on
$\mathbb{R}$ with $\otimes_{\mu}$ acting on $\mathcal{P}_s$ and
$\varkappa(\mu) >0$. A measure $\lambda\in\mathcal{P}_s$ is $\mu
$-weakly infinitely divisible if and only if there exists $A\geq0$ and
a symmetric $\sigma$-finite measure $\nu$ on $\mathbb{R}_0$ such
that $\nu ( [-a,a]^c  ) < \infty$ for each $a>0$,
\[
\int_0^{\infty} \mu \bigl( [-s, s]^c
\bigr) \nu(\mathrm{d}s) < \infty
\]
and
\[
\int_{\mathbb{R}}\mathrm{e}^{\mathrm{i}tx} (\mu\circ\lambda) (
\mathrm{d}x) = \exp \biggl\{- A|t|^{\varkappa(\mu)} - \int_{\mathbb{R}_0}
\bigl( 1 - \widehat {\mu}(ts) \bigr) \nu(\mathrm{d}s) \biggr\} .
\]
\end{theorem}

For details of the proof see \cite{weaklevy}. The parameter $A$ and
the measure $\nu$ we call the \textit{scale parameter} and $\mu
$-\textit{weak generalized L\'evy measure} respectively.
Bellow we present some examples of $\mu$-weakly stable distributions.
Since we consider symmetric measures, it is enough to restrict the
corresponding spectral measure $\nu$ to the positive half-line.

\subsection*{Examples}

{\renewcommand{\theexample}{3.4b}
\begin{example}\label{ex3.4b}
 Consider the Kendall weak generalized
convolution $\diamond_{\alpha}$, $\alpha\in(0,1]$, on $\mathcal
{P}_{+}$ defined by the measure $\mu_{\alpha} \in\mathcal{P}_s$
with the characteristic function $\widehat{\mu_{\alpha}}(t) = (1 -
|t|^{\alpha} )_{+}$ and the characteristic exponent $\varkappa(\mu
_{\alpha}) = \alpha$. We know (for details see, e.g., \cite{weaklevy}) that for every $p\leqslant\alpha$ there exists a
probability measure $\nu_{\alpha,p} \in\mathcal{P}_{+}$ such that
$\gamma_p = \mu_{\alpha} \circ\nu_{\alpha,p}$. The density of
$\nu_{\alpha,p}$ (which is ${\mu_{\alpha}}$-weakly $p$-stable) for
$p <\alpha$ is given by
\[
g_{\alpha,p}(s) = p \alpha^{-1} \bigl( (\alpha- p)
s^{-p-1} + p s^{-2p -1} \bigr) \mathrm{e}^{- s^{-p}}
\mathbf{1}_{(0,\infty)}(s).
\]
In the same paper \cite{weaklevy} it was shown that
\[
\exp \bigl\{ - |t|^p \bigr\} = \exp \biggl\{ - \int
_0^{\infty} \bigl( 1 - \bigl(1- |ts|^{\alpha}
\bigr)_{+} \bigr) \frac{p
(\alpha- p)}{\alpha s^{p+1}} \,\mathrm{d}s \biggr\}.
\]
Thus the L\'evy measure for symmetric $p$-stable measure with the
characteristic function $\exp \{ - |t|^p  \}$ can be
written as $\mu_{\alpha} \circ\lambda_{p}$, where $\lambda
_{p}(\mathrm{d}s)=p (\alpha- p)\alpha^{-1} s^{-p-1} \mathbf{1}_{(0,\infty
)}(s) \,\mathrm{d}s$. For $p=\alpha$ such a measure $\lambda_{\alpha}$ does not
exist, but we have that
\[
\exp \bigl\{ - |t|^{\alpha} \bigr\} = \lim_{p \nearrow\alpha} \exp
\biggl\{ - \int_0^{\infty} \bigl( 1 - \widehat{
\mu_{\alpha}} (ts) \bigr) \frac{p (\alpha- p)}{\alpha s^{p+1}} \,\mathrm{d}s \biggr\}.
\]
\end{example}}
{\renewcommand{\theexample}{3.5b}

\begin{example}\label{ex3.5b}
Consider the weakly stable Kingman distributions
\[
\omega_{s,1} (\mathrm{d}x) = \frac{\Gamma(s+1)}{\sqrt{\pi} \Gamma(s+ \sfrac
{1}{2})} \bigl( 1 -
x^2 \bigr)^{s-\sfrac{1}{2}} \mathbf{1}_{(-1,1)} (x) \,
\mathrm{d}x,
\]
$s > - \frac{1}{2}$, with the characteristic exponent $\varkappa
(\omega_{s,1}) = 2$. There exists a probability measure $\nu_{s,2}$
such that $\omega_{s,1} \circ\nu_{s,2} = N(0,1)$, where the density
of $\nu_{s,2}$ is given by
\[
f_{s,2}(x) = \frac{1}{2^s \Gamma(s+1)} x^{2s+1}
\mathrm{e}^{-{{x^2}/2}} \mathbf{1}_{(0,\infty)}(x).
\]
If by $\lambda_p$ we denote the distribution of the random variable
$\sqrt{\Theta}$, where $\Theta$ is the positive ${p/2}$-stable
random variable with the Laplace transform $\exp\{ - (2 t)^{ p/2} \}$,
then $\omega_{s,1} \circ\nu_{s,2} \circ\lambda_p = N(0,1) \circ
\lambda_p$ is symmetric $p$-stable. For $p<2$ the spectral measure for
$\gamma_p$ is a scale mixture of $\omega_{s,1}$ since for a suitable
constant $K>0$
\[
|t|^p = \int_0^{\infty} \bigl( 1 -
\widehat{\omega_{s,1}}(\operatorname{tr}) \bigr) \frac{K}{r^{p+1}} \,
\mathrm{d}r.
\]
\end{example}}

%s4 #&#
\section{L\'evy and additive processes with respect to
generalized and weak generalized convolutions}\label{sec4}

In this section, we consider an analog of a process with independent
increments, when the usual convolution is replaced by a generalized
one. To see that this is a natural generalization, consider the usual
process with independent increments $X=\{X_t: t\ge0\}$. $X$ is also a
Markov process with transition probabilities $P_{s,t}(x, \cdot)=
\delta_x \ast\lambda_{s,t}$, where probability measures $\lambda
_{s,t} = \mathcal{L}(X_t-X_s)$ satisfy an obvious consistency
condition: $\lambda_{s,t} \ast\lambda_{t,u}= \lambda_{s,u}$, $s<t <
u$. Conversely, given a family of distributions $\{\lambda_{s,t} \}$
satisfying the above consistency condition, there is a Markov process
$X$ with transition probabilities $P_{s,t}(x, \cdot)= \delta_x \ast
\lambda_{s,t}$. Due to the consistency condition, the increments of
$X$ are independent and determined by $\lambda_{s,t}$.
Therefore, the existence of a process with independent increments
follows from a standard construction of a Markov process with given
transition probabilities (see, e.g., Theorems 9.7 and 10.4 in
\cite{Sato}).

This approach was also applied by Nguyen Van Thu \cite{Thu} in the
context of generalized convolutions, and for Kingman's
convolutions in particular, to relate generalized L\'evy processes
to Bessel processes.

We will use this approach to define and construct additive processes
for generalized and weak generalized convolutions. We will identify
properties of convolutions that are needed for this construction to go
through, which indicates possible extensions beyond the types of
convolutions considered in this paper.
The consistency condition stated above naturally extends to the case of
generalized convolutions as follows
%
%e4.1 #&#
\begin{equation}
\label{cons} \lambda_{s,t} \diamond\lambda_{t,u}=
\lambda_{s,u}\qquad \forall s<t < u.
\end{equation}
It turns out that, given \eqref{cons} and properties of generalized
convolutions,
%
%e4.2 #&#
\begin{equation}
\label{C-K} P_{s,t}(x, \cdot) := \delta_x \diamond
\lambda_{s,t}(\cdot) ,\qquad s<t, x \in\mathbb{R}_+,
\end{equation}
satisfy the Chapman--Kolmogorov equations (see Theorem~\ref{Lp} below),
hence generalized additive process can be well-defined.

%
%de4.1 #&#
\begin{definition}
$X=\{X_t: t\ge0\}$ is said to be a $\diamond$-additive process
(associated with $\{\lambda_{s,t}\}$ satisfying
\eqref{cons}) if $X$ is a Markov process with transition
probabilities given by \eqref{C-K}. If $\lambda_{s,t}=\lambda
^{\diamond(t-s)}$ for some $\diamond$-infinitely
decomposable measure and all $0\le s<t$, then $X$ is called a $\diamond
$-L\'evy process generated by $\lambda$. The definition of $\otimes
_{\mu}$-additive and $\otimes_{\mu}$-L\'evy processes is analogous,
we replace $\diamond$ in the above by $\otimes_{\mu}$.
\end{definition}

The next theorem is stated in a greater generality to show that only
minimal conditions on convolutions are needed for the existence of
generalized additive processes.

%
%th4.2 #&#
\begin{theorem}\label{Lp}
Let $\mathbb{E}$ be a Polish space. Let $\star$ be a binary
associative operation on $\mathcal{P}(\mathbb{E})$ such that the map
$\mathbb{E}^{2} \ni(x,y) \mapsto\delta_x \star\delta_y(A) \in
[0,1]$ is measurable for each $A \in\mathcal{B}(\mathbb{E})$, and
for every $\lambda_1, \lambda_2 \in\mathcal{P}(\mathbb{E})$
%
%e4.3 #&#
\begin{equation}
\label{kernel3} \lambda_1 \star\lambda_2(A) = \int
_{\mathbb{E}^2} (\delta_x \star \delta_y) (A)
\lambda_1(\mathrm{d}x) \lambda_2(\mathrm{d}y).
\end{equation}
Given a family $\{\lambda_{s,t}: 0\le s<t \} \subset \mathcal
{P}(\mathbb{E})$ such that
\[
\lambda_{s,u} = \lambda_{s,t} \star\lambda_{t,u},
\qquad s<t<u,
\]
the probability kernels $P_{s,t}(x, \cdot) := \delta_x \star\lambda
_{s,t}(\cdot)$ on $\mathbb{E} \times\mathcal{B}$ satisfy the
Chapman--Kolmogorov equations, that is, for every $0<s<t<u$, $x\in
\mathbb{E}$ and $A\in\mathcal{B}(\mathbb{E})$,
%
%e4.4 #&#
\begin{equation}
P_{s,u}(x,A)=\int_{\mathbb{E}}P_{s,t}(x,
\mathrm{d}y)P_{t,u}(y,A) .\label{eq:Ch-K}
\end{equation}
Consequently, for any $\mu_0 \in\mathcal{P}(\mathbb{E})$,
there exists a Markov process $X=\{X_t \colon t \geqslant0\}$ in
$\mathbb{E}$ such that $\mathcal{L}(X_0)=\mu_0$ and, for all $t>s$,
$x \in\mathbb{E}$,
%
%e4.5 #&#
\begin{equation}
\label{cond} \mathbf{P} \bigl(X_t \in(\cdot) | X_s=x
\bigr) = \delta_x \star\lambda _{s,t}(\cdot).
\end{equation}
\end{theorem}

\begin{pf} Let $s, t, u, x$ and $A$ be as in \eqref{eq:Ch-K}. We have\vspace*{1pt}
\begin{eqnarray*}
P_{s,u}(x,A) &= &\delta_x \star \lambda_{s,u}(A)
= \delta_x \star (\lambda_{s,t} \star \lambda_{t,u})
(A)
\\
& =& (\delta_x \star \lambda_{s,t}) \star
\lambda_{t,u}(A)
\\
& =& \int_{\mathbb{E}^2} (\delta_y \star
\delta_z) (A) (\delta_x \star \lambda_{s,t}) (
\mathrm{d}y) \lambda_{t,u}(\mathrm{d}z)
\\
& =& \int_{\mathbb{E}^3} (\delta_x \star
\lambda_{s,t}) (\mathrm{d}y) (\delta _w \star
\delta_z) (A) \delta_y(\mathrm{d}w)
\lambda_{t,u}(\mathrm{d}z)
\\
& =& \int_{\mathbb{E}} (\delta_x \star
\lambda_{s,t}) (\mathrm{d}y) (\delta _y \star
\lambda_{t,u}) (A)
\\
& = &\int_{\mathbb{E}} P_{s,t}(x,\mathrm{d}y)P_{t,u}(y,A),
\end{eqnarray*}
where the third equality uses the associativity of $\star$; we also
applied \eqref{cons}--\eqref{kernel3}. The existence of the process
$X$ with desired properties follows now from \eqref{eq:Ch-K} by
Kolmogorov's extension theorem.
\end{pf}

%
%re4.3 #&#
\begin{remark}
Given a probability kernel $\mathbb{E}^2 \ni(x,y) \mapsto\rho_{x,y}
\in\mathcal{P}(\mathbb{E})$, one can define a ``convolution'' on
$\mathcal{P}(\mathbb{E})$ setting $\delta_x \star\delta_y :=\rho
_{x,y}$, and then extending $\star$ to arbitrary measures by \eqref{kernel3}.
If $\mathbb{E}$ is also a semigroup (not-necessarily commutative),
then it is natural to assume that $\rho_{x,0}=\rho_{0,x}=\delta_x$.
If $(\delta_x \star\delta_y) \star\delta_z = \delta_x \star
(\delta_y \star\delta_z)$ for all $x,y,z \in\mathbb{E}$, then
$\star$ is associative on $\mathcal{P}(\mathbb{E})$. In this way,
new classes of Markov processes, which are L\'evy processes relative to
such convolutions, can be defined.
\end{remark}

%
%th4.4 #&#
\begin{theorem}\label{th4.4}
Let $\star$ denote either a generalized convolution $\diamond$ or a
weak generalized convolution $\otimes_{\mu}$. Then for any consistent
family of probability measures $\{\lambda_{s,t}: 0\le s<t\}$ there
exists a $\star$-additive process $X=\{X_t: t\ge0\}$ generated by
this family and starting from 0.
If $\lim_{t \downarrow s} \lambda_{s,t} = \delta_0$ for every $s \ge
0$ [$\lim_{s \uparrow t} \lambda_{s,t} = \delta_0$ for every $t>0$,
resp.], then $X$ is right [left, resp.] continuous in probability.
Any $\star$-L\'evy process is continuous in probability.
\end{theorem}

\begin{pf} Suppose that $\lambda_{s,t} \to \delta_0$ as $t
\downarrow s$. For every $\epsilon>0$, by \eqref{cond} we have
\begin{eqnarray*}
\mathbf{P}\bigl( |X_t-X_s|> \epsilon\bigr) &=& \int\mathbf{P}\bigl(
|X_t-x|> \epsilon | X_s=x\bigr) \mathcal{L}(X_s)
(\mathrm{d}x)
\\
&=&\int\delta_x \star\lambda_{s,t} \bigl(\bigl\{y: |y-x| >
\epsilon\bigr\} \bigr) \mathcal{L}(X_s) (\mathrm{d}x)
\\
& \to&\int\delta_x \star\delta_0 \bigl(\bigl\{y: |y-x| >
\epsilon\bigr\} \bigr) \mathcal {L}(X_s) (\mathrm{d}x) = 0
\end{eqnarray*}
as $t \downarrow s$.
Similarly we treat continuity from the left. Now, if $X$ is a L\'evy
process, then the continuity of $\lambda_{s,t}= \lambda^{\star
(t-s)}$ follows from Proposition~\ref{prop1} and the beginning of
Section~\ref{sec3.2}.
\end{pf}

%
%re4.5 #&#
\begin{remark}\label{h-tail}
The $\diamond$-L\'evy processes are Markov processes in classical
sense. By Theorem~2.6 in \cite{Thu94}, it follows that if $\diamond$
is a generalized convolution on $\mathbb{R}_{+}$, or a weak
generalized convolution with $\varkappa(\diamond)>0$, then each
$\diamond$-L\'evy processes has strong Markov property, the Feller
property, it is continuous in probability and has c{\`a}dl{\`a}g
trajectories. Consequently, for each such process starting from a fixed
(nonrandom) point the Blumenthal's 0--1 law holds (see, e.g.,
Proposition~40.4 in \cite{Sato}).

Moreover, $\diamond$-L\'evy processes have heavy-tailed distributions
in each of the examples considered in this paper, provided $\varkappa
(\diamond) <2$ and $\diamond$ is not the maximum or stable
convolution. To see this it is enough to notice that in these cases for
all $x,y \in\mathbb{R} \setminus\{0\}$ the measure $\delta_x
\diamond\delta_y$ has infinite $p$-moment for $p > \varkappa
(\diamond)$. Such processes provide interesting new models for the
study of heavy-tail phenomena and possible long range dependence.
\end{remark}

%s5 #&#
\section{Stochastic integral processes with respect to \texorpdfstring{$\diamond$}{$diamond$}-L\'evy processes}\label{sec5}

For $\lambda$ being $\diamond$-infinitely decomposable probability
measure with the $\diamond$-generalized characteristic function
\[
\Phi_{\lambda} (t) = \exp \biggl\{ - A t^{\varkappa(\diamond)} - \int
_{0}^{\infty} \frac{1- \Omega(tx)}{ \upsilon(x) } m (\mathrm{d}x)
\biggr\},
\]
let $\mathcal{A}_{\lambda}$ be the class of nonnegative functions
$f$ on the positive half-line which are nonnegative, measurable,
bounded on compact intervals, and such that for every $t,u >0$
\[
\int_0^t f(x)^{\varkappa(\diamond)} \,\mathrm{d}x <
\infty,\qquad \int_0^{\infty} \int
_0^t \frac{1-\Omega(u f(x)s)}{\upsilon(s)} \,\mathrm{d}x m(
\mathrm{d}s) < \infty.
\]
By $\{X_t \colon t \geqslant0\}$ we denote the $\diamond$-additive
process based on $\lambda$ defined in the previous section.
We want to define a stochastic process
\[
Y_t = \diamond\int_0^t f(s) \,
\mathrm{d}X_s,\qquad t \geqslant0
\]
as a Markov process with the transition probabilities $P_{s,t}^f (x,
\cdot) = \delta_x \diamond P_{s,t}^f (0, \cdot)$ defined by the
$\diamond$-generalized characteristic function of $P_{s,t}^f (0, \cdot)$:
\begin{eqnarray*}
\Psi(f,s,t,u)
 =\exp \biggl\{ - A u^{\varkappa(\diamond)} \int_s^t
f(x)^{\varkappa(\diamond)} \,\mathrm{d}x - \int_0^{\infty}
\int_s^t \frac{1-
\Omega(u f(x)s) }{\upsilon(s)} \,\mathrm{d}x m(
\mathrm{d}s) \biggr\}.
\end{eqnarray*}
In view of the previous section the construction will be completed when
we prove the following:

%
%le5.1 #&#
\begin{lemma} For each $f \in\mathcal{A}_{\lambda}$ and every $s,t
\geqslant0$, $s<t$ the function $\Psi(f,s,t, \cdot)$ is a $\diamond
$-generalized characteristic function of a $\diamond$-infinitely
decomposable measure $P_{s,t}^f$.
\end{lemma}

\begin{pf}
Assume first that $f$ is a simple function, which
means that $f(x) = \sum_{k=1}^n a_k \mathbf{1}_{B_k}(x)$, where
$B_j\cap B_k = \emptyset$ for $j \neq k$ and $\bigcup_{k=1}^n B_k =
[s,t]$. We define the following measure
\[
P_{s,t}^f := T_{a_1} \lambda^{\diamond\ell(B_1)}
\diamond\cdots \diamond T_{a_n} \lambda^{\diamond\ell(B_n)},
\]
where $\ell$ is the Lebesgue measure on the positive half-line. We see that
\begin{eqnarray*}
\Phi_{P_{s,t}^f} (u) &=& \int_0^{\infty}
\Omega(ux) P_{s,t}^f (\mathrm{d}x)
\\
& =& \exp \Biggl\{ - A u^{\varkappa(\diamond)}\sum_{k=1}^n
a_k^{\varkappa(\diamond)} \ell(B_k) - \sum
_{k=1}^n \ell(B_k) \int
_{0}^{\infty} \frac{1- \Omega(ua_k x) }{ \upsilon(x) } m (\mathrm{d}x)
\Biggr\}
\\
& =& \exp \biggl\{ - A u^{\varkappa(\diamond)}\int_s^t
f(r)^{\varkappa(\diamond)} \,\mathrm{d}r - \int_s^t
\int_{0}^{\infty} \frac
{1- \Omega(u f(r) x) }{ \upsilon(x) } m (\mathrm{d}x)
\,\mathrm{d}r \biggr\}.
\end{eqnarray*}
Now, if $f \in\mathcal{A}_{\lambda}$, then there exists a sequence
of simple functions $f_n$ monotonically increasing to $f$ in each point
$r\in[s,t]$. By the Lebesgue dominated convergence theorem we have that
\begin{eqnarray*}
\lim_{n\rightarrow\infty} \Phi_{P_{s,t}^{f_n}} (u)
 =\exp \biggl\{ - A u^{\varkappa(\diamond)}\int_s^t
f(r)^{\varkappa
(\diamond)} \,\mathrm{d}r + \int_s^t
\int_{0}^{\infty} \frac{ \Omega(u f(r)
x) - 1 }{ \upsilon(x) } m (\mathrm{d}x)
\,\mathrm{d}r \biggr\}.
\end{eqnarray*}
Since the sequence of continuous functions converging to a continuous
function is converging uniformly on every compact interval, by the
definition of $\diamond$-generalized characteristic function there
exists a probability measure $P_{s,t}^f (0, \cdot)$ such that
$P_{s,t}^{f_n} \rightarrow P_{s,t}^f $ weakly if $n \rightarrow\infty
$ and
\[
\int_0^{\infty} \Omega(ux) P_{s,t}^f
( \mathrm{d}x) = \Psi(f,s,t,u).
\]
Infinite decomposability follows from the fact that
\[
\Psi(f,s,t,u) = \Psi_{A,m} (f,s,t,u) = \Psi_{{A/n}, {m/n}}^n
(f,s,t,u).
\]
It can also be derived from the following property:
\[
P_{s,t}^f \diamond P_{t,u}^f =
P_{s,u}^f,\qquad s<t<u.
\]
\upqed
\end{pf}

From this lemma, we conclude the following theorem.

%
%th5.2 #&#
\begin{theorem}
Let $\lambda$ be $\diamond$-infinitely decomposable and let $X=\{X_t
\colon t \geqslant0\}$ be the corresponding $\diamond$-L\'evy process
associated with $\lambda$. For given $f \in A_{\lambda}$, there
exists nonhomogenous Markov process $Y=\{Y_t \colon t \geqslant0\}$
with transition probabilities $\delta_x \diamond P_{s,t}^f$, where
$P_{s,t}$ are transition probabilities of $X$. The process $Y$ is a
$\diamond$-additive process which is denoted by
\[
Y_t = \diamond\int_0^t f(s) \,
\mathrm{d}X_s,\qquad t\geqslant0.
\]
\end{theorem}

%s6 #&#
\section{Weak generalized summation}\label{sec6}

Naturally, one would like to describe a generalized convolution in
terms of an operation on independent random variables. To this aim, one
can consider a weak generalized summation $X\oplus Y$ of nonnegative
random variables, where $\oplus$ is a binary operation on $\mathbb
{R}_+$. It turns out that this method is very restrictive, only
convolutions described in Examples \ref{ex2.3} and \ref{ex2.6} can be realized this
way. Indeed, if we assume that for all $a,b,c \geqslant0$, $a \oplus b
= b \oplus a$, $a \oplus0 = a$, $a \oplus(b \oplus c) = (a \oplus b)
\oplus c$ and $c (a \oplus b) = (ca) \oplus(cb)$, together with an
assumption on continuity, then by Bohnenblust's theorem (see \cite{bohn}), for some $\alpha\in(0,\infty]$,
\[
a \oplus b = \lleft\{ %
\begin{array} {l@{\qquad}l} \bigl(
a^{\alpha} + b^{\alpha} \bigr)^{1/{\alpha}} & \mbox{if } \alpha<
\infty,
\\\noalign{\vspace*{2pt}}
\max\{a,b\} & \mbox{if } \alpha= \infty. \end{array} %
\rright.
\]
%
%Considering the generalized adding on the set of all real random
%variables we could eventually extend this to the following
%$$
%a \oplus b = \left\{ \begin{array}{lcl}
% \left( a^{<\alpha>} + b^{<\alpha>} \right)^{<1/{\alpha}>} & \mbox{if}
%& \alpha<\infty, \\
% \max\{|a|,|b|\} {\rm sign}(a+b - {\rm min}(|a|, |b|) & \mbox{if} &
%\alpha= \infty,
%\end{array} \right.
%$$
%where $a^{<\alpha>} = |a|^{\alpha} {\rm sign}(a)$.

The problem of describing a generalized convolution in the language of
random variables seems to be difficult. However, weak generalized
convolutions open some new possibilities in this direction.

Recall that the random vector $\mathbf{X}$ and its distribution $\mu$
is weakly stable if for all random variables $\theta_1, \theta_2$ and
$\mathbf{X}_1, \mathbf{X}_2$ independent copies of $\mathbf{X}$ such
that $\theta_1, \theta_2, \mathbf{X}_1, \mathbf{X}_2$ are
independent there exists a random variable $\theta$ independent of
$\mathbf{X}$ such that
{\renewcommand{\theequation}{$\ast\ast$}
\begin{eqnarray}
\label{astast} \mathbf{X}_1 \theta_1 +
\mathbf{X}_2 \theta_2 \stackrel{\mathrm{d}} {=} \mathbf{X}
\theta.
\end{eqnarray}}
Until now, we were satisfied by defining the weak generalized
convolution based on this property:
\[
\mathcal{L}(\theta_1) \otimes_{\mu} \mathcal{L}(
\theta_2) = \mathcal{L}(\theta).
\]
Now we want to use the original property in defining weak generalized
addition which involves all random elements appearing in (\ref{astast}).

%
%le6.1 #&#
\begin{lemma}
Let $\mu$ be a nontrivial weakly stable distribution. Suppose that
$\mathbf{X}, \mathbf{X}_1$ and $\mathbf{X}_2$ are i.i.d. with
distribution $\mu$.
Then for all nonnegative random variables $\theta_1, \theta_2$ such that
$\theta_1, \theta_2, \mathbf{X}_1, \mathbf{X}_2$ are independent
there exist random elements $\mathcal{X}$, $\Theta$, $\mathcal{X}
\stackrel{d}{=} \mathbf{X}$, such that
\[
\theta_1 \mathbf{X}_1 + \theta_2
\mathbf{X}_2 = \mathcal{X} \cdot \Theta\qquad \mbox{a.e.}
\]
\end{lemma}

\begin{pf}
Let $\theta_1, \theta_2, \mathbf{X}_1, \mathbf
{X}_2$ be as assumed in the lemma, with random elements taking values
in a separable Banach space $\mathbb{E}$. By weak stability of
$\mathbf{X}$ we have that there exists independent random variable
$\Theta$ independent of $\mathbf{X}$ such that
\[
\theta_1 \mathbf{X}_1 + \theta_2
\mathbf{X}_2 \stackrel{\mathrm{d}} {=} \mathbf{X} \Theta.
\]
Corollary~5.11 in \cite{Kallenberg} states that for each two Borel
spaces $S$ and $T$, a measurable mapping $f\dvtx  T \rightarrow S$ and
some random elements $\xi$ in $S$ and $\eta$ in $T$ with $\xi
\stackrel{\mathrm{d}}{=} f(\eta)$ there exists a random element $\widetilde
{\eta}\stackrel{\mathrm{d}}{=} \eta$ in $T$ with $\xi= f(\widetilde{\eta
})$ a.e. We see that it is enough to apply this corollary for $\xi=
\theta_1 \mathbf{X}_1 + \theta_2 \mathbf{X}_2$, $\eta= (\mathbf
{X}, \Theta)$ and $f \dvtx \mathbb{E} \times[0,\infty) \mapsto
\mathbb{R}$ given by $f(\mathbf{x},s) = \mathbf{x}s$ to obtain
existence of $\widetilde{\eta} = (\mathcal{X}, \Theta)$ such that
\[
\theta_1 \mathbf{X}_1 + \theta_2
\mathbf{X}_2 = \mathcal{X} \Theta \qquad \mbox{a.e.}
\]
\upqed
\end{pf}

In the following definition by $\mathbb{K}$ we understand one of the
sets $\mathbb{R}$ or $\mathbb{R}_{+} =[0,\infty)$. If $\mu$ is
symmetric, then we can take $\mathcal{P}(\mathbb{K}) = \mathcal
{P}_s$ as well as $\mathcal{P} (\mathbb{K}) = \mathcal{P}_{+}$.

%
%de6.2 #&#
\begin{definition}\label{de6.2}
Let $(\Omega, \mathfrak{F}, \mathbf{P})$ be a rich enough
probability space, $\mu\in\mathcal{P}(\mathbb{E})$ be a nontrivial
weakly stable distribution, and let $s,t \in\mathbb{K}$.

The weak generalized convolution algebra $(\mathcal{P}(\mathbb{K}),
\otimes_{\mu})$ is representable (or the weak generalized convolution
$\otimes_{\mu}$ is representable) if there exist measurable functions
\[
\Theta\dvtx (\mathbb{K}\times\mathbb{E})^2 \rightarrow\mathbb{K}
\quad \mbox{and}\quad \mathcal{X} \dvtx (\mathbb{K}\times\mathbb{E})^2
\rightarrow\mathbb{E}
\]
such that for every choice of i.i.d. vectors $(\mathbf{X}_i)_{i \in
\mathbb{N}}$ with distribution $\mu$ and all $ i \neq j$, $i,j \in
\mathbb{N}$, the following conditions hold
\renewcommand\theenumi{(\arabic{enumi})}
\renewcommand\labelenumi{\theenumi}
\begin{enumerate}[(3)]
\item\label{q1} $\Theta(s, \mathbf{X}_i ; t, \mathbf{X}_j) =
\Theta(t, \mathbf{X}_j; s, \mathbf{X}_i)$ and $\mathcal{X}(s,
\mathbf{X}_i; t, \mathbf{X}_j) = \mathcal{X}(t, \mathbf{X}_j; s,
\mathbf{X}_i)$;
\item\label{q2} $\mathcal{X}(s, \mathbf{X}_i; t, \mathbf
{X}_j)\stackrel{d}{=} \mathbf{X}_1$;
\item\label{q3} $\mathcal{L} (\Theta(s, \mathbf{X}_i; t,
\mathbf{X}_j)  ) = \delta_s \otimes_{\mu} \delta_t$;
\item\label{q4} $ \Theta(s, \mathbf{X}_i; t, \mathbf{X}_j)$ and
$\mathcal{X}(s, \mathbf{X}_i; t, \mathbf{X}_j)$ are independent;
\item\label{q5} $s \mathbf{X}_i + t \mathbf{X}_j = \Theta(s,
\mathbf{X}_i; t, \mathbf{X}_j)\mathcal{X}(s, \mathbf{X}_i; t,
\mathbf{X}_j)$ a.e.;
\item\label{q6} $\Theta (\Theta(s, \mathbf{X}_i; t, \mathbf
{X}_j), \mathcal{X}(s, \mathbf{X}_i; t, \mathbf{X}_j); u, \mathbf
{X}_k )
= \Theta (s, \mathbf{X}_i; \Theta(t, \mathbf{X}_j; u, \mathbf
{X}_k), \mathcal{X}(t, \mathbf{X}_j; u, \mathbf{X}_k)  )$ a.e.
and
$\mathcal{X}  (\Theta(s, \mathbf{X}_i; t, \mathbf{X}_j),
\mathcal{X}(s, \mathbf{X}_i; t, \mathbf{X}_j); u, \mathbf{X}_k
 )
= \mathcal{X}  (s, \mathbf{X}_i; \Theta(t, \mathbf{X}_j; u,
\mathbf{X}_k), \mathcal{X}(t, \mathbf{X}_j;\allowbreak  u, \mathbf{X}_k)
)$ a.e.;
\item\label{q7} If $\sum_{i=1}^{\infty} s_i \mathbf{X}_i$
converges a.e. for some choice of $s_i \in\mathbb{K}$, $i \in\mathbb
{N}$, then $S_n \rightarrow S$ a.e. and $\mathcal{X}_n \rightarrow
\mathbf{X}$ a.e., where
\[
\begin{array} {c} S_1 = s_1,\qquad S_{n+1} =
\Theta( S_n, \mathcal{X}_n ; s_{n+1}, \mathbf
{X}_{n+1})\qquad \mbox{a.e.};
\\
\mathcal{X}_1 = \mathbf{X}_1, \qquad
\mathcal{X}_{n+1} = \mathcal{X}( S_n, \mathcal{X}_n
; s_{n+1}, \mathbf{X}_{n+1})\qquad \mbox{a.e.} \end{array}
\]
\end{enumerate}
\end{definition}

%
%le6.3 #&#
\begin{lemma} Assume that the weak generalized convolution algebra
$(\mathcal{P}(\mathbb{K}), \otimes_{\mu})$ is representable. If
$\theta_1, \theta_2$ are independent with distributions $\lambda_1,
\lambda_2$ respectively, and they are independent of $\mathbf{X}_1,
\mathbf{X}_2$, then the random elements $\Theta(\theta_1, \mathbf
{X}_1; \theta_2, \mathbf{X}_2)$ and $\mathcal{X}(\theta_1, \mathbf
{X}_1; \theta_2, \mathbf{X}_2)$ are independent.
\end{lemma}

\begin{pf} In fact, the result follows from properties (2) and (4)
of Definition~\ref{de6.2} by the following arguments:
\begin{eqnarray*}
&&\mathbf{P} \bigl\{ \mathcal{X}(\theta_1, \mathbf{X}_1;
\theta_2, \mathbf{X}_2) \in B, \Theta(
\theta_1, \mathbf{X}_1; \theta_2,
\mathbf{X}_2) \in A \bigr\}
\\
& &\quad = \int_{\mathbb{K}} \int_{\mathbb{K}}
\mathbf{P} \bigl\{ \mathcal{X}(s, \mathbf{X}_1; t,
\mathbf{X}_2) \in B, \Theta(s, \mathbf{X}_1; t,
\mathbf{X}_2) \in A \bigr\} \lambda_1(\mathrm{d}s)
\lambda _2 (\mathrm{d}t)
\\
& &\quad \stackrel{\mathrm{(4)}} {=} \int_{\mathbb{K}} \int
_{\mathbb{K}} \mathbf {P} \bigl\{\mathcal{X}(s, \mathbf{X}_1;
t, \mathbf{X}_2) \in B \bigr\} \mathbf{P} \bigl\{ \Theta(s,
\mathbf{X}_1; t, \mathbf {X}_2) \in A \bigr\}
\lambda_1(\mathrm{d}s) \lambda_2 (\mathrm{d}t)
\\
& &\quad \stackrel{\mathrm{(2)}} {=} \mu(B) \int_{\mathbb{K}} \int
_{\mathbb{K}} \mathbf{P} \bigl\{ \Theta(s, \mathbf{X}_1;
t, \mathbf{X}_2) \in A \bigr\} \lambda_1(\mathrm{d}s)
\lambda_2 (\mathrm{d}t)
\\
& & \quad = \mathbf{P} \bigl\{ \mathcal{X}(\theta_1,
\mathbf{X}_1; \theta_2, \mathbf{X}_2) \in B
\bigr\} \mathbf{P} \bigl\{ \Theta (\theta_1, \mathbf{X}_1;
\theta_2, \mathbf{X}_2) \in A \bigr\}.
\end{eqnarray*}
\upqed
\end{pf}

The following are examples of weak generalized convolutions that are
representable.

\subsection*{Examples}
{\renewcommand{\theexample}{6.1}
\begin{example}\label{ex6.1}
The symmetric convolution as the convolution on
$\mathcal{P}_{+}$ is representable and we have
\begin{eqnarray*}
\Theta\dvtx (\mathbb{R}_{+} \times \mathbb{R})^2 &
\rightarrow& \mathbb{R}_{+},\qquad \Theta(s,x;t,y) = | sx + ty |,
\\
\mathcal{X} \dvtx (\mathbb{R}_{+} \times\mathbb{R})^2 &
\rightarrow &\mathbb{R}_{+}, \qquad \mathcal{X} (s,x;t,y) =
\frac{ sx + ty}{| sx + ty|} = \sign(sx + ty).
\end{eqnarray*} %
\end{example}}

{\renewcommand{\theexample}{6.3}
\begin{example}\label{ex6.3} For $p\in(0,2]$ the weak generalized convolution
algebra $(\mathcal{P}_{+}, \ast_p)$ generated by symmetric
$p$-stable, weakly stable distribution $\gamma_p$ is evidently representable:
\begin{eqnarray*}
\Theta\dvtx (\mathbb{R}_{+} \times \mathbb{R})^2 &
\rightarrow& \mathbb{R}_{+}, \qquad \Theta(s,x;t,y) = \bigl\llVert (
s,t) \bigr\rrVert _p,
\\
\mathcal{X} \dvtx (\mathbb{R}_{+} \times\mathbb{R})^2 &
\rightarrow& \mathbb{R}_{+}, \qquad \mathcal{X} (s,x;t,y) =
\frac{s}{\llVert (
s,t)\rrVert _p} x + \frac{t}{\llVert ( s,t)\rrVert _p} y.
\end{eqnarray*} %
\end{example}}

{\renewcommand{\theexample}{6.5}
\begin{example}\label{ex6.5}
 The weak generalized convolution algebra
$(\mathcal{P}_{+}, \otimes_{\omega_n})$ is representable. The
corresponding functions are the following
\begin{eqnarray*}
\Theta\dvtx \bigl(\mathbb{R}_{+} \times \mathbb{R}^{n}
\bigr)^2 &\rightarrow& \mathbb{R}_{+},\qquad \Theta(s,
\mathbf{x}; t, \mathbf{y}) = \llVert s \mathbf{x} + t \mathbf{y}\rrVert
_2;
\\
\mathcal{X} \dvtx \bigl(\mathbb{R}_{+} \times\mathbb{R}^{n}
\bigr)^2 &\rightarrow&\mathbb{R},\qquad \mathcal{X} (s, \mathbf{x};
t, \mathbf{y}) = \frac{s \mathbf{x} + t \mathbf{y}}{\llVert  s \mathbf{x} + t
\mathbf{y}\rrVert _2}.
\end{eqnarray*} %
\end{example}}

If a nontrivial weak generalized convolution $\otimes_{\mu}$ is
representable and this will not lead to misunderstanding, we use the notation
\[
\Theta(\theta_1, \mathbf{ X}_1, \theta_2,
\mathbf{X}_2) = \theta _1 \oplus_{\mu}
\theta_2.
\]
In most of the cases, we shall however write
\[
\Theta(\theta_1, \mathbf{X}_1, \theta_2,
\mathbf{X}_2) = ( \theta_1 | \mathbf{X}_1 )
\oplus_{\mu} ( \theta_2 | \mathbf{X}_2 )
\]
and
\[
( \theta_1 | \mathbf{X}_1 ) \oplus_{\mu} (
\theta_2 | \mathbf{X}_2 ) \oplus_{\mu} \cdots
\oplus_{\mu} ( \theta_n | \mathbf{X}_n ) =: {\sum
_{i\leqslant
n}}^{\oplus_{\mu}} ( \theta_i |
\mathbf{X}_i ).
\]

To see the advantage of introducing representability for the weak
generalized convolution consider examples constructed as follows:

Let $\mathbf{X}$ with distribution $\mu$ be weakly stable and such
that the weak generalized convolution $\otimes_{\mu}$ is
representable. As in Section~\ref{sec4}, for any distribution $\lambda$ there
exists a Markov process $\{S_n \colon n \in\mathbb{N}_0\}$ with the
transition probabilities
\[
P_{n,k}(x, \cdot) = \delta_x \otimes_{\mu}
\lambda^{\otimes_{\mu
} (k-n)}.
\]
The existence of the process $\{S_n \colon n \in\mathbb{N}_0\}$
follows from a kind of existence theorem. Using the representability we
can do it more explicitly:

Let $\theta_i, i \in\mathbb{N}$, be a sequence of i.i.d. random
variables with distribution $\lambda$ and $\mathbf{X}_i, i \in
\mathbb{N}$ be a sequence of i.i.d. vectors with distribution $\mu$.
Now we define
\[
S_n := {\sum_{i\leqslant n}}^{\oplus_{\mu}} (
\theta_i | \mathbf{X}_i ), \qquad \mathbf{Z}_n
:= \sum_{i\leqslant n} \mathbf {X}_i
\theta_i = S_n \mathcal{X}_n,
\]
where $\mathcal{X}_1 = \mathbf{X}_1$, $\mathcal{X}_{n+1} = \mathcal
{X}( \theta_{n+1}, \mathbf{X}_{n+1}; S_n, \mathcal{X}_n)$. We see
that the sequence $\{S_n \mathcal{X}_n \colon n \in\mathbb{N}\}$ is
a classical independent increments homogenous random walk with the step
distribution $\mu\circ\lambda$. Considering simultaneously both
processes $\{ (S_n, \mathbf{Z}_n ) \colon n \in\mathbb{N}\}$ or even
all three processes $\{ (S_n, \mathcal{X}_n, \mathbf{Z}_n ) \colon n
\in\mathbb{N}\}$ we obtain more information than considering them separately.

\subsection*{Examples}

{\renewcommand{\theexample}{6.5a}
\begin{example}\label{ex6.5a}
In the case of $\mu= \omega_d$ uniform
distribution on the unit sphere $S_{d-1} \subset\mathbb{R}^d$ and
$\lambda$ with the density function
\[
f_{d-1,2}(r) = \frac{1}{2^{d-1} \Gamma(d)} r^{2d-1}
\mathrm{e}^{- {{r^2}/2}} \mathbf{1}_{(0,\infty)}(r)
\]
we see that $Z_n, n \in\mathbb{N}$, is the classical Wiener process
describing the position of the particle in $\mathbb{R}^d$ observed in
discrete times, $S_n, n \in\mathbb{N}$, describes the actual distance
of the particle form the origin, and the stationary process $\mathcal
{X}_n, n \in\mathbb{N}$ describes the projection of the actual
position of the particle on the unit sphere in $\mathbb{R}^d$.
\end{example}}

{\renewcommand{\theexample}{6.3a}
\begin{example}\label{ex6.3a}Another interesting example is connected with the
symmetric $\alpha$-stable L\'evy motion, where symmetry means in fact
spherical symmetry of the distribution of increments. To see this
notice first that every zero mean Gaussian random vector $\mathbf{X}$
is weakly stable and defines the representable weak generalized
convolution $\ast_2$:
\[
a \mathbf{X} + b \mathbf{X}' = \sqrt{a^2 +
b^2} \biggl( \frac
{a}{\sqrt{a^2 + b^2}} \mathbf{X} + \frac{b}{\sqrt{a^2 + b^2}}
\mathbf{X}' \biggr).
\]
Thus we have $\Theta(a, \mathbf{X}; b, \mathbf{X}') = \|(a,b)\|_2$, and
\[
\mathcal{X} \bigl(a, \mathbf{X}; b, \mathbf{X}' \bigr) = \biggl(
\frac{a}{\|
(a,b)\|_2} \mathbf{X} + \frac{b}{\|(a,b)\|_2} \mathbf{X}'
\biggr).
\]
We see that for the sequence $\mathbf{X}_n, n \in\mathbb{N}$, of
i.i.d. random vectors with rotationally invariant Gaussian distribution
and $\theta_i, i \in\mathbb{N}$, i.i.d. sequence of random variables
such that $\theta_i^2$ has $\frac{\alpha}{2}$-stable distribution
with the Laplace transform $\mathrm{e}^{-t^{\alpha/2}}$ the sequence $Z_n = S_n
\mathcal{X}_n$ consists of variables with symmetric $\alpha$-stable
distribution.
\end{example}}

Consequently the sequences $S_n, \mathcal{X}_n $ appearing in the
condition (7) of Definition~\ref{de6.2} are such that
\begin{eqnarray*}
S_n &:=& \Biggl( \sum_{i=1}^{n}
\theta_i^2 \Biggr)^{1/2} \qquad \mbox{is a
square root of a positive }\frac{\alpha}{2}\mbox{-stable process},
\\
Z_n &:=& \sum_{i\leqslant n} X_i
\theta_i = S_n \mathcal{X}_n \qquad \mbox{is
S$ \alpha$S rotationally invariant L\'evy process}.
\end{eqnarray*}

%
%le6.4 #&#
\begin{lemma}
The sequence $\sum_{1 \leqslant i\leqslant n}^{\oplus_{\mu}}  (
\theta_i | \mathbf{X}_i  ) $ converges a.e. if and only if the
sequence\linebreak[4]  $\sum_{1 \leqslant i\leqslant n} \theta_i \mathbf{X}_i$
converges a.e.
\end{lemma}

\begin{pf} Assume that the sequence $\sum_{1 \leqslant
i\leqslant n}^{\oplus_{\mu}}  ( \theta_i | \mathbf{X}_i
)$ converges a.e. (and in particular weakly) to a random variable
$\theta$. Since
\[
\biggl({\sum_{i\leqslant n}}^{\oplus_{\mu}} (
\theta_i | \mathbf{X}_i ) \biggr) \mathcal{X}_n
= \sum_{i=1}^n \theta _i
\mathbf{X}_i\qquad \mbox{a.e.}
\]
we see that the right-hand side of this equality converges weakly to a
random variable with distribution $\mathcal{L}(\theta) \circ\mu$.
Since the summands $\theta_i X_i$ are independent, the L\'evy's
equivalence theorem implies that $\sum_{1 \leqslant i\leqslant n}
\theta_i X_i$ converges a.e.

The opposite implication is a direct consequence of the property (7) of
representable weak generalized convolution.
\end{pf}

%s7 #&#
\section{Random measures with weak generalized summation}\label{sec7}

Let $(\mathbb{S},\mathcal{E})$ be a measurable space equipped with a
$\sigma$--finite measure $\varrho$. We define
\[
\mathcal{E}_0 = \bigl\{ A \in\mathcal{E} \colon\varrho(A) < \infty
\bigr\}.
\]
By $L^0(\Omega, \mathbb{E})$ we denote the space of all random
elements on $\Omega$ taking values in a separable Banach space
$\mathbb{E}$.

%
%de7.1 #&#
\begin{definition}
Let $\mu\in\mathcal{P}(\mathbb{E})$ be a nontrivial weakly stable
measure with representable convolution $\otimes_{\mu}$ and let
$\lambda\in\mathcal{P}$ be ${\mu}$-weakly infinitely divisible
measure. The set function
\[
\mathbf{M}_{\varrho,\lambda, \mu} = \mathbf{M} \dvtx \mathcal {E}_0
\rightarrow L^0(\Omega; \mathbb{R})\times L^0(\Omega;
\mathbb{E})
\]
is called the $\mu$-weak generalized random measure on a measurable
space $(\mathbb{S},\mathcal{E})$ with the control measure $\varrho$
if the following conditions hold:
\begin{enumerate}[(3)]
\item[(1)] $\mathbf{M}(\emptyset) = (0,\overline{ 0} )$ a.e.,
\item[(2)] $\mathbf{M}(A) =  ( \mathbf{M}_{\mu}(A), \mathbf
{M}_{\mu}(A) \mathcal{Y}(A)  )$, where $\mathbf{M}_{\mu}(A)$
has the distribution $\lambda^{\otimes_{\mu} \varrho(A)}$,
$\mathcal{Y}(A)$ has the distribution $\mu$, $\mathbf{M}_{\mu}(A)$
and $\mathcal{Y}(A)$ are independent for every set $A\in\mathcal{E}_0$,
\item[(3)] if the sets $A_1,A_2,\ldots,A_n\in\mathcal{E}_0$ are
disjoint, then the random vectors $\mathbf{M}(A_1), \mathbf{M}(A_2),\break
\ldots,  \mathbf{M}(A_n)$ are independent,
\item[(4)] if sets $A_1, A_2, \dots\in\mathcal{E}_0$ are disjoint
and $\bigcup_{i\in\mathbb{N}} A_i \in\mathcal{E}_0$, then
\[
\mathbf{M} \biggl( \bigcup_{i\in\mathbb{N}} A_i
\biggr){=} \biggl( \mathbf{M}_{\mu} \biggl( \bigcup
_{i\in\mathbb{N}} A_i \biggr), \sum
_{i \in\mathbb{N}}\mathbf{M}_{\mu}(A_i) \mathcal{Y}
( A_i) \biggr)\qquad \mbox{a.e.}
\]
\end{enumerate}
\end{definition}

For simplicity, we use the following notation
\[
\mathbf{M}_{\mu} \biggl( \bigcup_{i\in\mathbb{N}}
A_i \biggr){=} {\sum_{i \in\mathbb{N}}}^{\oplus_{\mu}}
\mathbf{M}_{\mu} (A_i ),
\]
when the sets $(A_i)$ are disjoint. Notice that this means that on the
second coordinate we have a random measure in the classical sense, that is,
the set function
\[
\bigl\{\mathbf{M}_{\ast}(A)= \mathbf{M}_{\mu}(A)
\mathcal{Y}(A) \colon A \in\mathcal{E}_0 \bigr\}
\]
is a classical independently scattered random measure and $\mathbf
{M}_{\ast}(A)$ has the distribution $(\mu\circ\lambda)^{\ast
\varrho(A)}$.

The existence of the $\mu$-weak generalized random measure can be
derived from the Kolmogorov extension theorem
by showing the consistency conditions for finite-dimensional
distributions of the process $\{\mathbf{M}(A):A\in\mathcal{E}_0\}$.
Instead of checking this directly we show that finite-dimensional
distributions of this process can be represented as distributions of
random vectors built with a collection of independent
two-dimensional random vectors and their $\oplus_{\mu}$-sums.

Let $A_1,\ldots,A_n\in\mathcal{E}_0$. Then there exist disjoint sets
$B_1,\ldots,B_N\in\mathcal{E}_0$ and
sets $I_1,\ldots,I_n\subseteq\{1,\ldots,N\}$ such that
\[
A_i=\bigcup_{j\in I_i}B_j,
\qquad i=1,\ldots,n.
\]
Consequently $\varrho(A_i) = \sum_{j\in I_i} \varrho(B_{j})$. For
each choice of $B_1, \dots, B_N \in\mathcal{E}_0$ we can choose
independent random variables $\theta_1, \dots, \theta_N$ with
distributions $\lambda^{\otimes_{\mu} \varrho(B_1)}, \dots,
\lambda^{\otimes_{\mu}\varrho(B_N)}$ respectively and a sequence of
i.i.d. random vectors $\mathbf{X}_1, \dots, \mathbf{X}_N$ with
distribution $\mu$ such that
\[
\Theta(A_i)= {\sum_{j\in I_i}}^{\oplus_{\mu}}
( \theta_j |\mathbf{X}_j )\qquad \mbox{a.e.}
\]
The random variables $\Theta(A_i)$ are well defined in view of
representability of the weak generalized convolution $\otimes_{\mu}$.
By the same argument, similarly as in the condition (6) in Definition~\ref{de6.2} we have uniquely, up to equality almost everywhere, defined vectors
$\mathbf{X}(A_i)$, $i=1,\dots, n$ such that $\Theta(A_i)$ and
$\mathbf{X}(A_i)$ are independent and
\[
\sum_{j\in I_i} \mathbf{X}_j
\theta_j = \mathbf{X}(A_i) \Theta (A_i)
\qquad \mbox{a.e.}
\]
Now it is easy to see that the random vector
\[
\bigl( \bigl(\Theta(A_1), \mathbf{X}(A_1)
\Theta(A_1) \bigr), \dots, \bigl(\Theta(A_n),
\mathbf{X}(A_n) \Theta(A_n) \bigr) \bigr)
\]
has the distribution desired for $(\mathbf{M}(A_1), \dots, \mathbf
{M}(A_n))$ and the consistency conditions in the Kolmogorov extension
theorem are evidently satisfied.

%s8 #&#
\section{L\'evy processes with respect to weak generalized summation}\label{sec8}

In this section, we assume that $\mathbb{S} = [0,\infty)$, $\mathcal
{E}= \mathcal{B}([0,\infty))$ and $\varrho$ is a $\sigma$-finite
measure on $[0,\infty)$, finite on compact sets.

%
%de8.1 #&#
\begin{definition}
Let $\mu$ be a weakly stable measure on $\mathbb{E}$ with
representable generalized convolution $\otimes_{\mu}$ and let
$\lambda\in\mathcal{P}$ be ${\mu}$-weakly infinitely
divisible. If $\mathbf{M}_{\varrho,\lambda, \mu}$
is the $\mu$-weak generalized random measure on $([0,\infty),\mathcal
{E})$ with the control measure $\varrho$, then the stochastic process
\[
\bigl\{ Z_{\varrho,\lambda, \mu} (t):= \mathbf{M}_{\mu} \bigl([0,t ) \bigr)
\colon t\geqslant0 \bigr\}
\]
is $\mu$-weakly additive that is, has $\mu$-weakly independent
increments, and the additive in classical sense process
\[
\bigl\{ \mathbf{Y}_{\varrho,\lambda, \mu} (t):= \mathbf{M}_{\mu} \bigl([0,t
) \bigr) \mathcal{Y} \bigl([0,t ) \bigr) \colon t\geqslant 0 \bigr\}
\]
is said to be associated with $  \{ Z_{\varrho,\lambda, \mu}
(t) \colon t\geqslant0  \}$.
\end{definition}

If the measures ${\varrho,\lambda, \mu}$ are fixed and this does not
cause a misunderstanding, then we use simplified notation
\[
\bigl\{ Z_{{\varrho,\lambda, \mu}}(t)\colon t \geq0 \bigr\} = \{ Z_t \colon
t \geq0 \}, \qquad \bigl\{ \mathbf{Y}_{{\varrho
,\lambda, \mu}}(t)\colon t \geq0 \bigr\} =
\{ \mathbf{Y}_t \colon t \geq0 \}.
\]
The finite dimensional distributions of $\mu$-weakly additive process
are uniquely determined by the measure $\varrho$ and the distribution
$\mathcal{L}(Z_1) = \lambda^{\otimes_{\mu} \varrho([0,1))}$.

%
%re8.2 #&#
\begin{remark}
Notice that if we have two independent $\mu$-weakly additive processes
$ \{ Z_{\varrho_1, \lambda_1, \mu}(t)\colon\allowbreak  t\geq0  \}
\mbox{ and }  \{ Z_{\varrho_2, \lambda_2, \mu}(t)\colon t\geq
0  \}$,
then their $\oplus_{\mu}$-sum
\[
\{ Z_t \colon t \geq0 \} = \bigl\{ Z_{\varrho_1,
\lambda_1, \mu}(t)
\oplus_{\mu} Z_{\varrho_2, \lambda_2, \mu}(t) \colon t\geq0 \bigr\}
\]
is also $\mu$-weakly additive in the following two cases:
\begin{enumerate}[(2)]
\item[(1)] if there exists a constant $a>0$ such that $\lambda
_1^{\otimes_{\mu}a} = \lambda_2$, and then
\[
Z_t = Z_{\varrho_1 + a \varrho_2, \lambda_1, \mu}(t),
\]
\item[(2)] if there exists $c>0$ such that $\varrho_2 = c\varrho_1$,
and then
\[
Z_t = Z_{\varrho_1, \lambda_1 \otimes_{\mu} \lambda_2^{\otimes
_{\mu}c}, \mu} (t).
\]
\end{enumerate}
\end{remark}

We want to consider a ${\mu}$-weakly additive process $ \{
Z_{\varrho, \lambda, \mu}(t) \colon t\geq0  \}$ as
a process with independent increments, but these increments shall not
be defined as a usual difference of
random variables. Thus, for every $0 \leq s\leq t$ we define the
increment between $Z_{s}$ and $Z_t$ to be
the random variable $Z_{s,t}= \mathbf{M}_{\mu}([s,t))$. By the
assumption, $Z_{s,t}$ is independent of $Z_s$ and
\[
Z_s \oplus_{\mu} Z_{s,t} = Z_t
\qquad \mbox{a.e.}
\]

%
%de8.3 #&#
\begin{definition}
A $\mu$-weakly additive stochastic process
\[
\bigl\{ Z_{\varrho,\lambda, \mu} (t) \colon t\geq0 \bigr\}
\]
is $\mu$-weak L\'evy process in law if the control measure $\varrho$
for the corresponding ${\mu}$-weak generalized random measure
$\mathbf{M}_{\varrho,\lambda,\mu}$ is proportional to the Lebesgue measure.
\end{definition}

It is easy to see that the stochastic process $\{\mathbf{Y}_t \colon t
\geq0\}$ associated with the $\mu$-weak L\'evy process $ \{
Z_{\varrho,\lambda, \mu} (t) \colon t\geq0  \}$ is a L\'evy
process in law in the classical sense.

A L\'evy process in law is an additive process with stationary
increments, that is continuous in probability. Since the control
measure $\varrho$ in the definition of the $\mu$-weak L\'evy process
in law is proportional to the Lebesgue measure, stationarity of
increments is evident. The next proposition implies that our process is
also continuous in probability.

%
%pr8.4 #&#
\begin{prop}
Let $\mu$ be a nontrivial weakly stable measure and let $\lambda$ be
$\mu$-weakly infinitely divisible. If the
measure $\varrho$ on $[0,\infty)$ does not have any atoms, then both
${\mu}$-weakly additive process $ \{
Z_{\varrho, \lambda, \mu}(t) \colon t\geq0  \}$ and the
process $ \{ \mathbf{Y}_t \colon t\geq0  \}$ associated
with $ \{ Z_{\varrho, \lambda,\mu}(t) \colon t\geq0  \}$
are continuous in probability.
\end{prop}

\begin{pf} Since $\varrho([s,t)) \rightarrow0$ for $t \searrow s$,
\[
\mathcal{L}(Z_{s,t}) = \lambda^{\otimes_{\mu} \varrho([s,t))} \rightarrow
\delta_0,
\]
which implies continuity in probability for the process $ \{
Z_{\varrho, \lambda, \mu}(t) \colon t\geq0  \}$.
Consequently, we have also
\[
\mathcal{L}(Y_t - Y_s) = \mu\circ\lambda^{\otimes_{\mu} \varrho
([s,t))}
\rightarrow\mu\circ\delta_0 = \delta_0 \qquad \mbox{for }
t \searrow s.
\]
\upqed
\end{pf}

%
%de8.5 #&#
\begin{definition}
Let $\mu\in\mathcal{P}(\mathbb{E})$ be a nontrivial weakly stable
measure and let $\ell$ be the Lebesgue measure on $[0,\infty)$. The
$\mu$-weak L\'evy process
\[
\bigl\{ N_{\mu}(t)\colon t \geq0 \bigr\} \stackrel{def} {=} \bigl\{
Z_{\ell, \lambda, \mu}(t)\colon t \geq0 \bigr\}
\]
is $\mu$-weak Poisson processes with the intensity $c>0$ if $\lambda=
\Exp_{\otimes_{\mu}}(c \delta_1)$.
\end{definition}

\subsection*{Examples}
{\renewcommand{\theexample}{8.3}
\begin{example}\label{ex8.3} Let $\mu= \gamma_p$, $p\in(0,2]$ be a
symmetric $p$-stable distribution on $\mathbb{R}$ with the
characteristic function $\mathrm{e}^{-A|r|^p}$, $A>0$. Then the $\mu$-weak
Poisson process $\{ N_{\gamma_p}(t)\colon t \geq0 \}$ is purely
discrete with the distribution
\[
\mathcal{L} \bigl(N_{\gamma_p}(t) \bigr) = \Exp_{\otimes_{\gamma_p}} (ct
\delta_1) = \mathrm{e}^{-ct} \sum
_{k=0}^{\infty} \frac{(ct)^k}{k!} \delta_{k^{1/p}}.
\]
It is easy to notice that the stochastic process $ \{ Y_t \colon
t\geq0  \}$ associated with $\{ N_{\gamma_p}(t)\colon t \geq0
\}$ is such that
\[
\Ee \mathrm{e}^{\mathrm{i}rY_t} = \widehat{\exp(ct \gamma_p)}(r) =
\exp \bigl\{ - ct \bigl( 1 - \widehat{\gamma_p}(r)\bigr) \bigr\} = \exp \bigl\{ -
ct \bigl( 1 - \mathrm{e}^{-A|r|^p} \bigr) \bigr\}.
\]
\end{example}}

{\renewcommand{\theexample}{8.4}
\begin{example}\label{ex8.4}
Consider the Kendall weak generalized convolution $\diamond_{\alpha}
\dvtx \mathcal{P}_s \rightarrow\mathcal{P}_s$ defined by the weakly
stable distribution
$\mu_{\alpha}$ on $\mathbb{R}$ with the characteristic function
$\widehat{\mu_{\alpha}}(t)=  ( 1-|t|^{\alpha}  )_+$,
$\alpha\in(0,1]$. By Example \ref{ex3.4as} we know that the
distribution of $N_{\mu_{\alpha}}(t)$ is given by
\[
\Exp_{\otimes_{\mu_{\alpha}}} (ct \delta_1 ) (\mathrm{d}s) =
\mathrm{e}^{-ct} ( \delta_0 + ct \lambda_0 ) (
\mathrm{d}s) + \frac{\alpha
(ct)^2}{2|s|^{ (2\alpha+1 )}} \mathrm{e}^{- ct|s|^{-\alpha}} \mathbf{1}_{(1,\infty)}\bigl(|s|\bigr)
\,\mathrm{d}s.
\]
The L\'evy stochastic process in law $ \{ Y_t \colon t\geq0
 \}$ associated with the $\mu_{\alpha}$-weak Poisson process
$ \{ N_{\mu_{\alpha}} (t) \colon t\geq0  \}$ is such that
\[
\mathbf{E} \mathrm{e}^{\mathrm{i}rY_t} = \exp \bigl\{ - ct \bigl( 1 - \widehat{
\mu _{\alpha}}(r) \bigr) \bigr\}= \mathrm{e}^{- ct |r|^{\alpha}}\mathbf
{1}_{[-1,1]}(r) + \mathrm{e}^{- ct }\mathbf{1}_{[-1,1]^c}(r).
\]
This means that
\[
\mathcal{L}(Y_t) (\mathrm{d}s) = \mathrm{e}^{- ct }
\delta_0 (\mathrm{d}s) + \bigl(1- \mathrm{e}^{-ct} \bigr)
f_{\alpha} (s) \,\mathrm{d}s,
\]
where
\[
f_{\alpha} (s) = \frac{1}{\pi} \frac{1}{1 - \mathrm{e}^{-ct}} \int
_0^{1} \cos(sr) \bigl( \mathrm{e}^{-ct |r|^{\alpha}}
- \mathrm{e}^{-ct} \bigr) \,\mathrm{d}r.
\]
For $\alpha= 1$ we obtain
\[
f_1(s) = \frac{ct}{\pi(\mathrm{e}^{ct}-1)} \frac{s\mathrm{e}^{ct} - ct \sin(s) - s
\cos(s)}{s((ct)^2 +s^2)}.
\]
\end{example}}

{\renewcommand{\theexample}{8.5}
\begin{example}\label{ex8.5} Consider $\{ N_{\omega_{3,1}}(t)\colon t \geq0
\}$, the $\omega_{3,1}$-weak Poisson process with the intensity $c>0$.
In this construction, we assume that $\otimes_{\omega_{3,1}} \colon
\mathcal{P}_s \rightarrow\mathcal{P}_s$.

The distribution of $N_{\omega_{3,1}}(t)$ we obtain substituting $c$
by $ct$ in the formula obtained in Example \ref{ex3.5a}, thus
\[
\mathcal{L} \bigl(N_{\omega_{3,1}}(t) \bigr)= \exp ( ct \omega_{3,1}
) \ast ( \delta_0-ct \omega_{3,1} + ct
\lambda_{0} ).
\]
The L\'evy process in law $\{Y_t \colon t \geq0\}$ associated with $\{
N_{\omega_{3,1}}(t)\colon t \geq0 \}$ is such that
\[
\mathcal{L}(Y_t) = \exp(ct \omega_{3,1}) =
\mathrm{e}^{-ct} \sum_{n=
0}^{\infty}
\frac{(ct)^n}{n!} \omega_{3,1}^{\ast n}.
\]
Since $\omega_{3,1}$ is the uniform distribution on $[-1,1]$, $\omega
_{3,1}^{\ast n}$ are also well known and, for example, in \cite{Kill} we can
find that $\omega_{3,1}^{\ast n}$ has the following density function:
\[
f^{(n)}(x) = \lleft\{ %
\begin{array} {l@{\qquad}l} { \displaystyle \sum
_{i=0}^k} (-1 )^i \displaystyle {{n}
\choose{i}}\displaystyle  \frac{ (x+n-2i )^{n-1}}{ (n-1
)! 2^n}, & x \in[ -n + 2k,-n + 2 (k+1 )),
\\
& k=0, \ldots, n-1,
\\
0, & \mbox{otherwise}.
\end{array} %
\rright.
\]
\end{example}}

%
%re8.6 #&#
\begin{remark}
The idea of a stochastic process associated with another process
suggests a
natural connections with the idea of subordinated processes described
in Feller's monograph \cite{Feller}. The construction there was
the following: We start with two independent stochastic processes $\{
X_t \in\mathbb{R} \colon t\geq0\}$ and $\{
\mathbb{T}(t)\in[0,\infty) \colon t\geq0\}$, $\mathbb{T}(t)$
increasing, and we define
\[
\{ X_{\mathbb{T}(t)} \colon t \geq0 \}.
\]
The process $\{ X_{\mathbb{T}(t)} \colon t\geq0\}$ is subordinated to
the process $\{X_t \colon t\geq0\}$ by $\{\mathbb{T}(t)
\colon t\geq0\}$. This construction is rich enough to cover many cases.

One of the best known subordinated processes is the sub-stable
independent increments process. It is based on a
strictly stable process $\{ X_t \colon t\geq0\}$ with independent
stationary increments. This means that
\[
X_t \stackrel{\mathrm{d}} {=} t^{1/{\alpha}} X_1,
\qquad X_{t+s} - X_t \stackrel {\mathrm{d}} {=}
s^{1/{\alpha}} X_1,\qquad X_t \bot(X_{t+s}
- X_t).
\]
The corresponding time stochastic process $\{ \mathbb{T}(t) \colon
t\geq0\}$ takes values in the positive half-line, has
independent increments and the Laplace transform $\mathbf{E} \mathrm{e}^{-r
\mathbb{T}_t} = \exp\{ - \operatorname{tr}^{\beta}\}$ for some $\beta<1$. Then
$ \{X_{\mathbb{T}(t)} \colon t \geq0 \}$ is an $(\alpha
\beta)$-stable stochastic process with independent increments.

The same process can be obtained by our construction as associated with
the $\mu$-weakly additive process
\[
\bigl\{ Z_{\varrho, \lambda, \gamma_{\alpha}}(t) = \mathbb {T}(t)^{1/{\alpha}}\colon t \geq0 \bigr
\},
\]
where
\[
\lambda= \mathcal{L} \bigl( \mathbb{T}(1)^{1/{\alpha}} \bigr),\qquad \mu=
\gamma_{\alpha} = \mathcal{L}(X_1),\qquad \varrho= \ell.
\]
In this case $Z_{s,t} = M_{\mu}([s,t)) = (\mathbb{T}(t) - \mathbb
{T}(s))^{1/{\alpha}}$ and $Y_t - Y_s \stackrel{\mathrm{d}}{=} Z_{s,t} X_1$.
Thus the associated process $\{ Y_t \colon t \geqslant0\}$ was
obtained by some operation on the space, not by randomizing the time,
as $\{ X_{\mathbb{T}(t)} \colon t \geq0\}$; however they are
stochastically equivalent. In the case $\alpha=2$ and $\{X_t \colon t
\geqslant0\}$ being multidimensional Brownian motion we again obtain
rotationally invariant independent increment symmetric $2\beta$-stable
stochastic process.
\end{remark}

%s9 #&#
\section{Weak stochastic integrals}\label{sec9}

In this section, we give a construction of a stochastic integral using
the weak generalized summation.
We assume that the considered nontrivial weakly stable measure $\mu$
belongs to $\mathcal{P}$ (instead of $\mu\in\mathcal{P}(\mathbb
{E})$ for the sake of simplicity) and that the weak generalized
convolution $\otimes_{\mu}$ is representable. Let $\lambda$ be $\mu
$-weakly infinitely divisible, $\mathbf{M}_{\varrho,\lambda, \mu}$
be $\mu$-weak generalized random measure for some $\sigma$-finite
measure $\varrho$ on $(\mathbb{S}, \mathcal{E})$, and let $\mathcal
{E}_0 = \{ A \in\mathcal{E}\colon\varrho(A)<\infty\}$.

The representability property of $\otimes_{\mu}$ allows us to
construct a stochastic integral as in the case of the usual convolution
(see, e.g., Rajput and Rosi\'nski \cite{RR}). We will only outline
this construction. For a simple function
\[
f(x) = \sum_{i=1}^n a_i
\mathbf{1}_{A_i} (x),
\]
where $A_1, \dots, A_n \in\mathcal{E}_0$ are disjoint sets and
$a_1,\dots,a_n \in\mathbb{R}$, put
\[
I_{\varrho, \lambda,\mu} (f) = \int_{\mathbb{S}} f(x) \mathbf
{M}_{\mu}(\mathrm{d}x) \stackrel{\mathrm{def}} {=} { \sum
_{i \leq n}}^{\oplus_{\mu}} a_i \mathbf{M}_{\mu}(A_i).
\]

%
%le9.1 #&#
\begin{lemma}
Assume that we have two representations for the simple function $f$, that is,
\[
f(x) = \sum_{i=1}^n a_i
\mathbf{1}_{A_i} (x)\quad \mbox{and}\quad f(x) = \sum
_{i=1}^m b_i \mathbf{1}_{B_i}
(x),
\]
such that $A_1, \dots, A_n, B_1, \dots, B_m \in\mathcal{E}_0$ and
$A_i \cap A_j = \emptyset$, $B_i \cap B_j = \emptyset$ for $i \neq
j$. Then
\[
{ \sum_{i \leq n}}^{\oplus_{\mu}} a_i
\mathbf{M}_{\mu}(A_i) = { \sum
_{j \leq m}}^{\oplus_{\mu}} b_j
\mathbf{M}_{\mu}(B_j)\qquad \mbox{a.e.}
\]
\end{lemma}

\begin{pf} There exists a family of disjoint sets $C_1,\dots,
C_N \in\mathcal{E}_0$ such that for every $i \leq n$ and $j \leq m$
there exists $I_i = \{k_{1,i}, \dots, k_{n_i, i}\} \subset\{1, \dots
, N\}$ and $J_j = \{ \ell_{1,j},\dots, \ell_{m_j,j} \} \subset\{1,
\dots, N\}$ such that
\[
\bigcup_{k \in I_i} C_k = A_i,
\qquad \bigcup_{\ell\in J_j} C_{\ell} =
B_j,\qquad i \leq n, j \leq m.
\]
Of course $I_k \cap I_i = \emptyset$ and $J_k \cap J_i = \emptyset$
for $k \neq i$. Thanks to representability of $\otimes_{\mu}$ it
makes sense to consider generalized sums, thus by our construction
\[
\mathbf{M}_{\mu} (A_i) = {\sum
_{k \in I_i}}^{\oplus_{\mu}} \mathbf {M}_{\mu}(C_k)
\qquad \mbox{a.e.}\quad \mbox{and}\quad \mathbf{M}_{\mu}
(B_j) = {\sum_{\ell\in J_j}}^{\oplus_{\mu}}
\mathbf{M}_{\mu}(C_{\ell}) \qquad \mbox{a.e.}
\]
Put $c_k := a_i = b_j$ if $C_k \subset A_i \cap B_j$.
Now we see that the following equalities hold almost everywhere
%
%e9.1 #&#
\begin{eqnarray*}
{ \sum_{i \leq n}}^{\oplus_{\mu}} a_i
\mathbf{M}_{\mu}(A_i) & = & { \sum
_{i \leq n}}^{\oplus_{\mu}} a_i {\sum
_{k \in I_i}}^{\oplus
_{\mu}} \mathbf{M}_{\mu}(C_k)
= { \sum_{i \leq n}}^{\oplus_{\mu}} {\sum
_{k \in I_i}}^{\oplus_{\mu}} a_i \mathbf{M}_{\mu}(C_k)
\nonumber
\\
& = & { \sum_{i \leq n}}^{\oplus_{\mu}} {\sum
_{k \in I_i}}^{\oplus
_{\mu}} c_k \mathbf{M}_{\mu}(C_k)
= { \sum_{k \leq N}}^{\oplus
_{\mu}} c_k
\mathbf{M}_{\mu}(C_k)
\nonumber
\\
& = & { \sum_{j \leq m}}^{\oplus_{\mu}} {\sum
_{\ell\in
J_j}}^{\oplus_{\mu}} c_{\ell}\mathbf{M}_{\mu}(C_{\ell})=
{ \sum_{j \leq m}}^{\oplus_{\mu}} b_j {
\sum_{\ell\in J_j}}^{\oplus_{\mu
}} \mathbf{M}_{\mu}(C_{\ell})
\nonumber
\\
& = & { \sum_{j \leq m}}^{\oplus_{\mu}}
b_j \mathbf{M}_{\mu}(B_j).
\end{eqnarray*}
\upqed
\end{pf}

%
%re9.2 #&#
\begin{remark}
Let $\mu= \mathcal{L}(\mathbf{X}) \in\mathcal{P}(\mathbb{E})$ be
a nontrivial weakly stable measure and let $\lambda= \Exp_{\otimes
_{\mu}} ( \delta_1)$. The $\mu$-weak generalized random measure
$\mathbf{M}_{\mu}$ consists of the variables with $\mu$-weak Poisson
distribution and $\mathcal{L}(\mathbf{M}_{\mu}(A)) = \Exp_{\otimes
_{\mu}} (\varrho(A) \delta_1)$ for $A \in\mathcal{E}_0$. Since the
$\otimes_{\mu}$-generalized characteristic function of the measure
$\Exp_{\otimes_{\mu}} ( a \delta_1)$ is equal to the classical
characteristic function of $\Exp_{\otimes_{\mu}} ( a \delta_1)
\circ\mu= \exp( a \delta_1 \circ\mu) = \exp(a \mu)$ then
\begin{eqnarray*}
\Ee\exp \bigl\{ \mathrm{i} \bigl\langle \mathbf{t}, a \mathbf{M}_{\mu}(A)
\mathbf{X}\bigr\rangle \bigr\} & = & \exp \bigl\{ - \bigl( 1 - \widehat{\mu} (a
\mathbf{t}) \bigr) \varrho(A) \bigr\}
\\
& = & \exp \biggl\{ - \int_{\mathbb{S}} \bigl( 1 - \widehat{\mu}
\bigl(a \mathbf{1}_A(x) \mathbf{t} \bigr) \bigr) \varrho(
\mathrm{d}x) \biggr\}.
\end{eqnarray*}
If $f(x)=\sum_{i=1}^{n}a_i\mathbf{1}_{A_i}(x)$,
for disjoint $A_1,\ldots,A_n\in\mathcal{E}_0$, then
\begin{eqnarray*}
\Ee\exp \bigl\{\mathrm{i}\bigl\langle \mathbf{t}, I_{\varrho,\lambda, \mu}(f) \mathbf
{X}\bigr\rangle \bigr\} &= &\prod_{k=1}^n
\Ee\exp \bigl\{\mathrm{i}\bigl\langle\mathbf{t}, a_k
\mathbf{M}_{\mu} (A_k )\mathbf{X}\bigr\rangle \bigr\}
\\
& = & \exp \Biggl\{ - \sum_{k=1}^n \int
_{\mathbb{S}} \bigl( 1 - \widehat{\mu} \bigl(a_k
\mathbf{1}_{A_k}(x) \mathbf{t} \bigr) \bigr) \varrho (\mathrm{d}x)
\Biggr\}
\\
& = & \exp \biggl\{ - \int_{\mathbb{S}} \bigl( 1 - \widehat{\mu}
\bigl( f(x) \mathbf{t} \bigr) \bigr) \varrho(\mathrm{d}x) \biggr\}
\\
& = & \exp \biggl\{ - \int_{\mathbb{S}} \bigl( 1 - \widehat{\mu} (
s \mathbf{t}) \bigr) \varrho_f(\mathrm{d}s) \biggr\},
\end{eqnarray*}
where $\varrho_f(A)=\varrho (f^{-1}(A) )=\varrho (\{
x\in\mathbf{E}:f(x)\in A\} )$, $A\in\mathcal{E}_0$. This
means that\vspace*{1pt}
\[
\mathcal{L} \bigl( I_{\varrho, \lambda, \mu}(f) \bigr) = \Exp _{\otimes_{\mu}} (
\varrho_f).
\]
\end{remark}

%
%pr9.3 #&#
\begin{prop}\label{pr9.3} Assume that the weakly stable measure $\mu= \mathcal
{L}(X)$ on $\mathbb{R}$ is nontrivial and symmetric with the
characteristic exponent $\varkappa$. Let $\lambda$ be $\mu$-weakly
infinitely divisible with the scale parameter $A\geq0$ and the $\mu
$-weak generalized L\'evy measure $\nu$. Let $f\dvtx  \mathbb{S} \mapsto
\mathbb{R}$ be a measurable function such that\vspace*{1pt}
\[
\int_{\mathbb{S}} \bigl|f(x)\bigr|^{\varkappa} \varrho(\mathrm{d}x) <
\infty \quad \mbox{and}\quad \int_{\mathbb{R}} \int
_{\mathbb{S}} \bigl\llvert 1 - \widehat{\mu } \bigl(f(x) ts \bigr)
\bigr\rrvert \varrho(\mathrm{d}x) \nu(\mathrm{d}s) < \infty.
\]
Then the stochastic integral $I_{\varrho, \lambda,\mu} (f)$ exists
as the limit in probability of stochastic integrals of simple functions.
Moreover, the $\otimes_{\mu}$-generalized characteristic function of
$I_{\varrho, \lambda,\mu} (f)$ is of the form\vspace*{1pt}
\begin{eqnarray*}
&&\Ee\exp \bigl\{ \mathrm{i} t I_{\varrho, \lambda,\mu} (f) X \bigr\}
\\
&&\quad = \exp \biggl\{ - A |t|^{\varkappa} \int_{\mathbb{S}}
\bigl|f(x)\bigr|^{\varkappa} \varrho(\mathrm{d}x) - \int_{\mathbb{R}} \int
_{\mathbb
{S}} \bigl( 1 - \widehat{\mu} \bigl(f(x) ts \bigr) \bigr)
\varrho(\mathrm{d}x) \nu (\mathrm{d}s) \biggr\}.
\end{eqnarray*}
\end{prop}

\begin{pf} It is enough to prove this for simple function $f=
\sum_{i=1}^n a_i \mathbf{1}_{A_i}$ for disjoint sets $A_1, \dots,
A_n$. Notice that the generalized characteristic function for $\lambda
$ is the following\vspace*{1pt}
\[
\widehat{\lambda\circ\mu}(t)=\exp \biggl\{-A|t|^{\varkappa(\mu
)}-\int
_{\mathbb{R}} \bigl(1-\widehat{\mu}(ts) \bigr)\nu(\mathrm{d}s) \biggr
\}.
\]
Since\vspace*{1pt}
\begin{eqnarray*}
\bigl(T_{a_1}\lambda^{\varrho(A_1)}\otimes_{\mu}\cdots
\otimes _{\mu} T_{a_n} \lambda^{\varrho(A_n)} \bigr)\circ
\mu&=& T_{a_1}\lambda^{\varrho
(A_1)}\circ\mu\ast\cdots\ast
T_{a_n} \lambda^{\varrho(A_n)}\circ \mu
\\
& =& T_{a_1} (\lambda\circ\mu )^{\ast\varrho(A_1)} \ast\cdots\ast
T_{a_n} (\lambda\circ\mu )^{\ast\varrho(A_n)},
\end{eqnarray*}
where, for the simplicity, we write $\lambda^{\varrho(A_i)}$ instead
of $\lambda^{\otimes_\mu\varrho(A_i)}$, we have\vspace*{1pt}
\begin{eqnarray*}
&&\Ee\exp \bigl\{ \mathrm{i} t I_{\varrho, \lambda,\mu} (f) X \bigr\}
\\
&&\quad =\prod_{i=1}^{n}\exp \biggl
\{-A|ta_i|^{\varkappa(\mu)}\varrho (A_i)-
\varrho(A_i)\int_{\mathbb{R}} \bigl(1-\widehat{
\mu}(a_its) \bigr)\nu(\mathrm{d}s) \biggr\}
\\
&&\quad = \exp \Biggl\{ - A \sum_{i=1}^{n}|ta_i|^{\varkappa}
\varrho(A_i) - \int_{\mathbb{R}} \sum
_{i=1}^n \bigl( 1 - \widehat{\mu}
(a_its)\varrho(A_i) \bigr) \nu(\mathrm{d}s) \Biggr\}
\\
&&\quad = \exp \biggl\{ - A |t|^{\varkappa} \int_{\mathbb{S}}
\bigl|f(x)\bigr|^{\varkappa} \varrho(\mathrm{d}x) - \int_{\mathbb{R}} \int
_{\mathbb
{S}} \bigl( 1 - \widehat{\mu} \bigl(f(x) ts \bigr) \bigr)
\varrho(\mathrm{d}x) \nu (\mathrm{d}s) \biggr\}.
\end{eqnarray*}
This ends the proof.
\end{pf}

%
%re9.4 #&#
\begin{remark}
The Proposition~\ref{pr9.3} states in particular that the random variable
$I_{\varrho, \lambda, \mu}(f)$ is $\mu$-weakly infinitely divisible
with the scale parameter
\[
A' = A \int_{\mathbb{S}} \bigl|f(x)\bigr|^{\varkappa(\mu)}
\varrho(\mathrm{d}x),
\]
and the $\mu$-weak generalized L\'evy measure $\varrho_f \circ\nu$,
where for $A \in\mathcal{E}_0$ $\varrho_f(A) = \varrho(f^{-1}(A))$.
\end{remark}

%\begin{appendix}
%\section{}
%\end{appendix}

% zodis "Acknowledgments" paliekamas pagal autoriu
\section*{Acknowledgements}
Jan Rosi\'nski's research was partially supported by a grant \# 281440
from the
Simons Foundation.

The authors are grateful to the anonymous referee for careful reading
of the
manuscript and constructive comments.

%\begin{supplement}%[id=suppA]
%\sname{Supplement A}
%\stitle{}
%\slink[doi]{10.3150/00-BEJXXXXSUPP} %[doi,text={...}] - jei reikia
%suskaldyti doi
%\sdatatype{.pdf}
%\sfilename{BEJ000\_supp.pdf}
%\sdescription{}
%\end{supplement}

% imsref loaded by jurgita.kaciuliene, 2014-07-30 13:32:04

\printhistory
\end{document}